\documentclass[a4paper,11pt]{article}

\usepackage[latin1]{inputenc}

\usepackage{epsfig}
\usepackage{amsmath,amsfonts,amssymb}
\usepackage{color}
\usepackage[hypertex]{hyperref}

\usepackage{mathrsfs} 

\input{macro}
\newcommand{\h}{\hbar}
\newcommand{\E}{\mathscr{E}}
\newcommand{\D}{\mathscr{D}}
\renewcommand{\O}{\mathscr{O}}
\renewcommand{\H}{\mathscr{H}}
\newcommand{\BH}{\textbf{H}}
\newcommand{\BPi}{\boldsymbol{\Pi}}
\renewcommand{\P}{\mathscr{P}}

\newcommand{\M}{\mathscr{M}}
\newcommand{\N}{\mathscr{N}}
\renewcommand{\B}{\mathscr{B}}
\newcommand{\R}{\mathcal{R}}
\newcommand{\p}{\mathfrak{p}}
\renewcommand{\P}{\mathcal{P}}

\newcommand{\bracket}[1]{\langle{#1}\rangle}
\newcommand{\sfrac}[2]{{\textstyle\frac{#1}{#2}}}

\newcommand{\Sp}{\operatorname{Sp}}  
\newcommand{\Id}{\operatorname{Id}}  

\title{Spectral asymptotics\\ via the semiclassical Birkhoff normal
  form}

\author{Laurent \textsc{Charles}\footnote{Institut de
    Math{\'e}matiques de Jussieu (UMR 7586), Universit{\'e} Pierre et
    Marie Curie -- Paris 6, Paris, F-75005 France.} \and \textsc{V\~u
    Ng\d oc} San\footnote{Institut Fourier (UMR 5582), Universit\'{e} Joseph
    Fourier -- Grenoble 1, BP 74, 38402-Saint Martin d'H\`{e}res
    Cedex, France.}}

\begin{document}
\maketitle

\begin{footnotesize}
  \noindent \textbf{MS Classification :}
  58J40, 
  58J50, 
  58K50, 
  47B35, 
  53D20, 
  81S10. 
\end{footnotesize}

\begin{abstract}
  This article gives a simple treatment of the quantum Birkhoff normal
  form for semiclassical pseudo-differential operators with smooth
  coefficients. The normal form is applied to describe the discrete
  spectrum in a generalised non-degenerate potential well, yielding
  uniform estimates in the energy $E$. This permits a detailed study
  of the spectrum in various asymptotic regions of the parameters
  $(E,\h)$, and gives improvements and new proofs for many of the
  results in the field.  In the completely resonant case we show that
  the pseudo-differential operator can be reduced to a Toeplitz
  operator on a reduced symplectic orbifold.  Using this quantum
  reduction, new spectral asymptotics concerning the fine structure of
  eigenvalue clusters are proved. In the case of polynomial
  differential operators, a combinatorial trace formula is obtained.
\end{abstract}

\section{Introduction}

The Birkhoff normal form, in classical mechanics, is a well known
refinement of the averaging method~: under a suitable canonical
transformation, a perturbation of a harmonic oscillator $H_2$ can be
replaced by its average along the classical Hamiltonian flow generated
by $H_2$. With the averaging method, this remains valid as long as one
restricts the dynamics to times bounded by $\O(1/\epsilon)$, where
$\epsilon$ is the size of the perturbation. Using the Birkhoff normal
form, this time can be extended to $\O(1/\epsilon^N)$ for arbitrary
$N$, provided one takes into account higher order terms which are also
averaged, but in a more intricate sense. Note that, in this work, we
do not try to impose special restrictions to the original Hamiltonian
that would imply some better convergence properties (Gevrey
convergence, or even analyticity). Instead, we take any smooth
function and perform the Birkhoff normal form in a neighbourhood of a
non-degenerate minimum.

In quantum mechanics, it is known since at least 1975 that an analogue
of the Birkhoff normal form can be applied in a very successful way.
At the formal level, this is attested by physicists
like~\cite{eckhardt,ali}. Adding on top of this the experience of
excellent numerical computations, it has become an important tool for
molecular physics~(see~\cite{swimm-delos} and more
recently~\cite{joyeux-birkhoff,robert-joyeux}).

On the mathematics side, the Birkhoff normal form for
pseudo-differential operators near a non-degenerate minimum of the
symbol has been used by several authors already. In particular the
article~\cite{sjostrand-semi} by Sj{\"o}strand is very interesting with
this respect, but only deals with the \emph{non-resonant} normal form,
that it, when the harmonic oscillator is of the form
\begin{equation}
  H_2(x,\xi)=\sum_{j=1}^n \nu_j(x_j^2+\xi_j^2)/2,
  \label{equ:intro-H2}
\end{equation}
where the coefficients $\nu_j$ are linearly independent over the
rationals. The result is that, when the energy $E$ is of order
$\h^\gamma$ with $\gamma\in(0,1)$ then, for $\h$ small enough, the
quantum system has the same spectrum as a completely integrable
Hamiltonian.

\ouf

The initial goal of our work here is to extend this to the resonant
case. But since our methods also give new proofs for Sj{\"o}strand's
result and some improvements, and moreover unify them with the
analysis of low-lying eigenvalues initially discovered by
Simon~\cite{simon1} and
Helffer-Sj{\"o}strand~\cite{helffer-sjostrand}, it might be of
interest to present it here in the general case.  Moreover, we believe
that several intermediate statements are of independent value, and
involve for the main part only standard results of semiclassical
analysis (symbolic and functional calculus for pseudo-differential
operators). In particular we are able to compare the initial
pseudo-differential operator to a differential operator with
polynomial coefficients, which is very important for many practical
purposes, including numerical computations. On the other hand, the
treatment of the resonant case is very hard to perform within the
standard pseudo-differential calculus because of unavoidable
singularities due to the fact that, when the coefficients $\nu_j$ are
integers, the $S^1$-action generated by the time-$2\pi$ flow of $H_2$
is in general non free~: periodic orbits with smaller periods appear.
This explains why so little results were available in that case. Our
strategy here is to abandon pseudo-differential operators for Toeplitz
operators, in the spirit of Boutet de Monvel and Guillemin~\cite{BG}.
The appropriate theory that can deal with orbifold singularities was
developed in~\cite{charles-toeplitz} and~\cite{charles-reduc}.

\ouf

Let us briefly describe our spectral result in this case. Assume
$P=-\frac{h^2}{2}\Delta+V(x)$ is a Schr{\"o}dinger operator with a smooth
potential $V$ on $X=\RM^n$ or on an $n$-dimensional compact manifold
$X$ equipped with a smooth density. (More generally, $P$ could be any
pseudo-differential operator in some standard class, which is actually
our assumption in this article.)  Assume $V\in\Cinf(X)$ has a global
minimum at a point $0$ which we shall call here the origin, and
suppose this minimum is non-degenerate. By a linear, unitary change of
variable in local coordinates near $0$, one can always assume that
$V''(0)$ is diagonal; let $(\nu_1^2,\dots,\nu_n^2)$ be its
eigenvalues, with $\nu_j>0$. The rescaling $x_j\mapsto
\sqrt{\nu_j}x_j$ transforms $P$ into a perturbation of the harmonic
oscillator $\hat{H}_2$:
\[
P=\hat{H}_2+W(x), \qquad \text{ with }\; \hat{H}_2=\sum_{i=1}^n
\frac{\nu_j}2\left(-\h^2\frac{\partial^2}{\partial x_j^2} +
  x_j^2\right),
\]
where $W(x)$ is a smooth potential of order $\O( |x|^3)$ at the origin.

Now assume that the coefficients $\nu_j$ are \emph{completely
  resonant}: there exist a real number $\nu_c>0$ and coprime positive
integers $\p_1,\dots,\p_n$ such that $\nu_j=\nu_c\p_j$. Then the
spectrum of $\hat{H}_2$ consists of the arithmetic progression
$E_N=\h\nu_c(\frac{\abs{\nu}}2+N)$ for $N\in\NM$, with multiplicity of
order $N^{n-1}$ as $N\fleche\infty$. It is then expected that, for
small energies, the spectrum of $P$ is a perturbation of the spectrum
of $\hat{H}_2$, splitting each eigenvalue $E_N$ into a band, or
\emph{cluster}. We prove this in a precise way. Actually, we describe
in theorem~\ref{theo:int1} the size and the internal structure of each
cluster, as follows. Let $H_2$ be the corresponding classical harmonic
oscillator, as in~(\ref{equ:intro-H2}). It has a $2\pi$ periodic
Hamiltonian flow $\phy_t$.  Let $k=k(x,\xi)$ be the average of $W$
along this flow.
\[
k(x,\xi) = \frac{1}{2\pi}\int_0^{2\pi}W\circ\varpi(\phy_t(x,\xi))
\]
where $\varpi:T^*X\fleche X$ is the cotangent projection. Let
$S_N\subset T^*X$ be the sphere:
\[
S_N=\{(x,\xi)\in T^*X,\quad H_2(x,\xi)=E_N \}.
\]

\vspace{1ex}
\noindent\textbf{Theorem }[theorem~\ref{theo:int1}]\emph{%
  \begin{enumerate}
  \item There exists $\hbar_0>0$ and $C>0$ such that for every $\hbar
    \in (0, \h_0]$
    \begin{equation*}
      \Sp ( P ) \cap (-\infty, C \hbar ^{\frac{2}{3}} ) \subset \bigcup_{E_N
        \in \Sp ( \hat{H}_2 )} \Bigl[ E_N - \frac{\nu_c \hbar}{3}, E_N +
      \frac{\nu_c\hbar}{3} \Bigr].
    \end{equation*}
  \item When $E_N\leq C\h^{\frac23}$, let $m (E_N, \h) = \# \Sp (P)
    \cap \Bigl[ E_N - \frac{\nu_c \hbar}{3}, E_N +
    \frac{\nu_c\hbar}{3} \Bigr]$. Then $m(E_N,\h)$ is precisely the
    dimension of $\ker(\hat{H}_2-E_N)$.
  \item Let $E_N+\lambda_1(E_N,\h),\dots,
    E_N+\lambda_{m(E_N,\h)}(E_N,\h)$ be the eigenvalues of $P$ in this
    $N$-eth band. Then, uniformly for $\h<\h_0$ and $N$ such that
    $E_N\leq C\h^{\frac23}$,
    \begin{align}
      \lambda_1 (E_N, \h) &= \inf_{(x,\xi) \in S_N}
      \abs{k(x,\xi)} + (E_N)^{\frac{3}{2}} \O ( N^{-1} ),\\
      \lambda_{m(E_N,\h)} (E_N, \h) &= \sup_{(x,\xi) \in S_N}
      \abs{k(x,\xi)} + (E_N)^{\frac{3}{2}} \O (N^{-1})
    \end{align}
    and for any function $g \in \Cinf ( \RM)$,
    \[
    \sum_{i=1}^{m (E_N, \h)} g \Bigl( \frac{\lambda_i (E_N, \h)}
    {(E_N)^{\frac{3}{2}}} \Bigr) = \Bigl( \frac{1}{2\pi
      \h}\Bigr)^{n-1} \!\!\!\int_{S_N} \!\!\!g \Bigl(
    \frac{k(x,\xi)}{(E_N)^{\frac{3}{2}}} \Bigr) \mu_{E_N} (x,\xi) +
    \O(N^{2-n})
    \]
    where $\mu_{E_N}$ is the Liouville measure of $S_N$
  \end{enumerate}}

Thus we see that the average perturbation $k$ behaves as a principal
symbol for the spectral analysis restricted to each cluster. 
Several
improvements of this statement are proved in the article. First, $k$
can actually be replaced by the homogeneous term of degree 3 in its
Taylor expansion. Secondly, the exponent $2/3$ in the term $C\h^{2/3}$
(and in $(E_N)^{3/2}$ where its inverse appears) is due to the fact
that in general resonances of order 3 may happen in $H_2$~: relations
of the form $\p_j=2\p_i$ or $\p_i=\p_j+\p_k$. If one rules these out,
then the exponent $2/3$ can be replaced by $1/2$, but in general with
a modified $k$ (if the homogeneous term of degree 3 in the Taylor
expansion of the potential vanishes, then $k$ keeps the same
definition. Otherwise the formula is more involved). Finally, the last
expansion in the theorem is actually the leading order of a full
asymptotic expansion in $\h/E$.  In particular, in sub-principal
terms, one can exhibit oscillatory contributions of type $\zeta^N$
where $\zeta$ is some (finite order) complex root of 1.

The estimate of the spectral density in the particular case $\p_1 =
...=\p_n =1$ was first obtained in the thesis of the second author
\cite{san-these} through a reduction to Toeplitz
operators. Independently Bambusi and Tagliaferro conjectured and
partially proved the estimates of the smallest and largest eigenvalues
in each band. Then Bambusi and the first author worked on a proof
using the quantum Birkhoff normal form of \cite{BGP}  and Toeplitz
operators. The result was announced in \cite{bambusi-talk}.

\ouf

On technical side, it might be worth mentioning here that we do not
use any exotic pseudo-differential calculus in order to deal with
formal Taylor series. Instead of this, we rely extensively on various
scaling properties of the harmonic oscillator $H_2$, which seems
particularly fit for this purpose. This allows us to play all along
with $(E,\h)$ as two (almost) independent small parameters. The
results of Sj{\"o}strand are thus recovered in the regime $E=\h^\gamma$.
Then, when we study the resonant case, these scaling properties of
$H_2$ become even more crucial, because the effective semiclassical
parameter becomes $h=\h/E$ (sections~\ref{sec:joint}
and~\ref{sec:toeplitz}).

\ouf

To conclude this introduction, let us mention that quantum Birkhoff
normal forms have also become a very important tool in inverse
spectral problems. Formally, the Birkhoff normal form is a
(semi)classical invariant from which, generally under analyticity
assumptions, one can hope to recover the full classical dynamics (see
for instance~\cite{guillemin-wave} and~\cite{zelditch-revolution}).
This aspect is not discussed here.

\paragraph{Structure of the article. ---}
The Birkhoff normal form is based on a simple formal construction,
which can be explained directly in a quantum setting; that's what we
recall in section~\ref{sec:formal} (theorem~\ref{theo:formal}). Most
of the material in this section is not new; however it is crucial
here, and the notation introduced there is used throughout the
article. As we next show in section~\ref{sec:semi}, the relevance of
the formal result to semiclassical operators is due to a general
theorem allowing to compare the spectrum of pseudo-differential
operators in the so-called semi-excited regime on the basis of the
Taylor expansions of the symbols (theorem~\ref{theo:taylor}). Adding
standard arguments of spectral theory we obtain a general statement of
the quantum Birkhoff normal form (theorem~\ref{theo:FNB}). \emph{Up to
  an error of size $\O(E^\infty)+\O(\h^\infty)$, it reduces the
  spectral problem to the analysis of a pseudo-differential operator
  $K$ commuting with a quantum harmonic oscillator $\hat{H}_2$.} In
section~\ref{sec:joint} we describe the joint spectrum of $P$ and its
Birkhoff normalisation $K$, proving an important estimate relating the
formal order of $K$ with its operator norm, when restricted to
eigenspaces of $\hat{H}_2$ (lemma~\ref{lemm:normes}). Several
applications are given~: Weyl asymptotics, expansions of the low-lying
eigenvalues, and the use of polynomial differential operators (thus
giving a rigorous justification of the spectroscopy computations
of~\cite{swimm-delos,joyeux-birkhoff}).  Finally the last
section~\ref{sec:toeplitz} is devoted to the resonant case, when
$\nu_j$ are integers, up to a common multiple. Then the classical flow
of $H_2$ is periodic, and it is known that one expects the spectrum to
exhibit clustering. We describe these clusters of eigenvalues.
Technically and conceptually, the main result is that the restriction
of $K$ to eigenspaces of $\hat{H}_2$ can be identified to a Toeplitz
operator on the corresponding reduced symplectic orbifold
(theorem~\ref{theo:toeplitz}). This allows us to introduce $\h/E$ as a
\emph{second semiclassical parameter} and yields spectral asymptotics
for these clusters in terms of the principal symbol of $K$
(theorem~~\ref{theo:int1}). A more precise trace formula involving
sub-orbifolds and hence oscillatory terms is given in~\ref{theo:int2}.
We end the article with an amusing combinatorial formula expressing a
certain sum over integral points of a rational polytope, which comes
as a direct consequence of our results (theorem~\ref{theo:polytope}).

\ouf

\paragraph{Acknowledgements ---} Laurent Charles thanks Dario
  Bambusi for collaborating on the subject.

\section{The formal Birkhoff normal form}
\label{sec:formal}

The Weyl quantisation on $\RM^{2n}=T^*\RM^n$ is based on a particular
grading for formal symbols in $x,\xi,\h$, where the degree in the
semiclassical parameter counts \emph{twice} the degree of each other
variable $x_i$ or $\xi_i$. This grading is particularly adapted to the
harmonic oscillator and hence to the quantum Birkhoff normal form. It
also appears naturally in the context of deformation
quantisation~\cite{fedosov-book}. We mainly follow here the
presentation of~\cite{san-these}, but other authors have used this
approach.

Thus we work with the space
\[
\E = \CM\formel{x_1,\dots,x_n,\xi_1,\dots,\xi_n,\h},
\]
and we define the weight of the monomial $x^\alpha\xi^\beta\h^\ell$ to
be $\abs{\alpha}+\abs{\beta}+2\ell$. The finite dimensional vector
space spanned by monomials of weight $N$ shall be denoted by $\D_N$.
Let $\O_N$ be the subspace consisting of formal series whose
coefficients of weight $<N$ vanish. $(\O_N)_{N\in\NM}$ is a filtration
\[
\E=\O_0 \supset \O_1\supset \cdots , \qquad \bigcap_N \O_N = \{0\},
\]
and shall be used for all convergences in this section.

The bracket associated to the Weyl product on $\E$ defines a Poisson
algebra structure on $\E$: it is the unique bilinear bracket for which
$\h$ is central, which satisfies the Jacobi identity, the Leibniz
identity (with the Weyl product), and which is commutative on all
generators amongst $x_1,\dots,x_n$, $\xi_1,\dots,\xi_n$ and $\h$,
except for the relations
\[
\forall j=1,\dots,n, \quad [\xi_j,x_j]=\frac{\h}{i}.
\]
Notice that this structure is invariant by linear canonical changes of
coordinates.  There is a simple formula (Moyal's formula) for the
brackets of two elements of $\E$, but we shall not need it in this
article. However in the following sections we will use the fact that
if $\hat{H}$ and $\hat{P}$ are Weyl-quantisations of symbols $H$ and
$P$ with formal Taylor series at the origin $[H]$ and $[P]$ in $\E$, then the Taylor
series of the symbol of the operator commutator $[\hat{H},\hat{P}]$ is
precisely the Weyl bracket $[[H],[P]]$.

The filtration of $\E$ has a nice behaviour with respect to the Weyl
bracket.  If $N_1+N_2\geq 2$ then
\[
\h^{-1}[\O_{N_1},\O_{N_2}] \subset \O_{N_1+N_2-2}.
\]
If $A\in\E$ the adjoint operator $P\mapsto[A,P]$ will be denoted by
$\ad{A}$. We shall be interested in the adjoint action of elements of
$\D_2$. Such elements are of the form $\h H_0 + H_2$, where
$H_0\in\CM$ and $H_2$ is a quadratic form in $(x,\xi)$. Since $\h$ is
central, we may restrict here to quadratic forms only.  They will be
called \emph{elliptic} when the quadratic form is positive. Because of
this grading we see that when $H_2\in\D_2$, then $\h^{-1}\ad{H_2}$
acts as an \emph{endomorphism} of each $\D_N$.  A fundamental property
of the Weyl bracket is that $\frac{i}{\h}\ad{H_2}P$ is exactly the
classical Poisson bracket $\{H_2,P\}$.

We will say that $H_2\in\D_2$ is \emph{admissible}
when $\D_N=\ker(\ad{H_2})+\textup{im}(\ad{H_2})$. A typical example is
the harmonic oscillator (see lemma~\ref{lemm:commutation} below):
\[
H_2 = \nu_1(x_1^2+\xi_1^2)/2+\cdots+\nu_n(x_n^2+\xi_n^2)/2.
\]
One can show that all elliptic $H_2$ can be written as harmonic
oscillators in some canonical coordinates and hence
are admissible as well. Indeed, eigenvalues of Hamiltonian matrices
come by pairs $(\nu_i,-\nu_i)$~: this implies that an elliptic $H_2$
must have the form of a harmonic oscillator plus some nilpotent terms.
But no such nilpotent term is allowed to show up because the flow of
$H_2$ is contained in the hypersurface $\{H_2=\textup{const}\}$, which
is compact.

The formal quantum Birkhoff normal form can be expressed as follows.
\begin{theo}
  \label{theo:formal}
  Let $H_2\in\D_2$ be admissible and $L\in\O_3$. Then there exists
  $A\in\O_3$ and $K\in\O_3$ such that
  \begin{itemize}
  \item $e^{i\h^{-1}\ad{A}}(H_2+L) = H_2+K$ ;
  \item $[K,H_2]=0$ .
  \end{itemize}
  Moreover if $H_2$ and $L$ have real coefficients then $A$ and $K$
  can be chosen to have real coefficients as well.
\end{theo}
Notice that the sum
\[
e^{i\h^{-1}\ad{A}}(H_2+L) = \sum_\ell
\frac{1}{\ell!}\left(\frac{i}{\h}\ad{A}\right)^\ell(H_2+L)
\]
is indeed convergent in $\E$ because $\frac{i}{\h}\ad{A}$ sends $\O_N$
into $\O_{N+1}$.

\begin{demo}
  We construct $A$ (and hence $K$) by successive approximations with
  respect to the filtration of $\E$. Modulo $\O_3$ the equality is
  trivial. So let $N\geq 1$ and suppose that for some $A_N\in\O_3$ we
  have
  \[
  e^{i\h^{-1}\ad{A_N}}(H_2+L) = H_2+K_3+\cdots+K_{N+1} + R_{N+2} +
  \O_{N+3},
  \]
  where $K_i\in\D_i$ and commutes with $H_2$, $R_{N+2}\in\D_{N+2}$.
  Let $A'\in\D_{N+2}$; then a small calculation gives
  \[
  e^{i\h^{-1}\ad{A_N+A'}}(H_2+L) = H_2+K_3+\cdots+K_{N+1} + K_{N+2} +
  \O_{N+3},
  \]
  where
  
  \begin{equation}
    K_{N+2}=R_{N+2}+i\h^{-1}\ad{A'}H_2 = R_{N+2} - i\h^{-1}\ad{H_2}A' .
    \label{equ:decomp}
  \end{equation}
  We look for an $A'$ such that $[K_{N+2},H_2]=0$. This is possible
  because $H_2$ is admissible. Now if we assume that $H_2$, $L$ and
  $K_j$, $j\leq N+1$ are real, then $R_{N+2}$ is real too. Since
  $\frac{i}{\h}\ad{H_2}=\{H_2,\cdot\}$ is a real endomorphism, we have
  \[
  \D_N^\RM = \ker^\RM(\ad{H_2})+\textup{im}^\RM(\ad{H_2})~.
  \]
  Hence\eqref{equ:decomp} can be solved with real coefficients.
\end{demo}

\begin{rema}
  \label{rema:class_birkhoff}
  If we write the theorem modulo $\h$ we recover the classical
  Birkhoff normal form for Hamiltonians on $\RM^{2n}$. Indeed, let $p$
  and $a$ be $\Cinf$ functions on $\RM^{2n}$ with Taylor expansion at
  the origin $H_2 (x, \xi) + L(x,\xi,0)$ and $A(x, \xi,0)$
  respectively. Then if $\phi$ denotes the Hamiltonian flow of $a$ at
  time $1$, we have
  $$ p \circ \phi = H_2 + k 
  $$
  where $k$ has the asymptotic expansion $K(x, \xi, 0)$. Consequently
  the Poisson bracket of $H_2$ and $k$ is flat at the origin.
\end{rema}

\begin{rema} 
  Another way of constructing the quantum Birkhoff normal form would
  be to start from the classical result and build successively in
  increasing powers of $\h$. This was used by several authors and
  amounts to follow a different filtration which, in a sense, is less
  optimal than ours.  Nevertheless the result is the same, as for
  instance in~\cite{sjostrand-semi}.
\end{rema}

\begin{rema}
  The result presented here is often called the Birkhoff-Gustavson
  normal form in the mathematical physics literature. Gustavson
  popularised the idea of Birkhoff in~\cite{gustavson} by providing
  computer programs performing the canonical transformation involved.
  Moreover, Gustavson added the analysis of the resonant cases, while
  in his treatise~\cite{birkhoff}, Birkhoff only dealt with the
  non-resonant situation. Note that Moser had a similar result before
  Gustavson, in the article~\cite{moser-new}. In some sense, the
  Birkhoff normal form is the Hamiltonian version of the
  Poincar{\'e}-Dulac method~\cite{dulac}.  Actually the Poincar{\'e}-Dulac
  normal form is even more general since it allows for the hypothesis
  of admissibility to be relaxed.  Then $H_2$ has to be split into
  commuting semisimple and nilpotent parts, and the normalisation is
  performed with respect to the semisimple part.  The quantum version
  is probably much more complicated to analyse, but it would be very
  interesting to do so.
\end{rema}

In this article we will always assume that $H_2$ is elliptic~: in some
canonical coordinates, one can write
\[
H_2=\sum_{j=1}^n \frac{\nu_j}{2}(x_j^2+\xi_j^2) .
\]

In order to understand what kind of formal series $K$ can show up in
the Birkhoff normal form, it is crucial to study the kernel of
$\ad{H_2}$. The following lemma is elementary and standard.
\begin{lemm}
  \label{lemm:commutation}
  $\ad{H_2}$ is diagonal on the $\CM\formel{h}$-basis
  $z^\beta\bar{z}^\gamma$ where $\beta,\gamma$ are multi-indices in
  $\NM^n$ and $z_j=x_j+i\xi_j$, and
  \begin{equation}
    \ad{H_2}(z^\beta\bar{z}^\gamma)=\pscal{\beta-\gamma}{\nu}
    z^\beta\bar{z}^\gamma
  \end{equation}
\end{lemm}
We also state and prove the following --- maybe less standard ---
result, which will be one of the tools in the next sections to obtain
a pseudo-differential version of the Birkhoff normal form.
Let $\R$ be the resonance module
  \begin{equation}
    \R:=\{\alpha\in\ZM^n,\quad \pscal{\alpha}{\nu}=0\}
    \label{equ:module}
  \end{equation}
and denote by $n-k$ its rank ($k \geq 1$). 

\begin{lemm}
  \label{lemm:torus}
  There exists a Hamiltonian $\T^k$ action on $\RM^{2n}=T^*\RM^n$
  such that the space of all power series that
  commute with $H_2$ is exactly the space of $\T^k$-invariant power
  series:
  \begin{equation}
    \ker \ad{H_2} = \E^{\T^k} .
  \end{equation}
\end{lemm}
\begin{demo}
  Let us use the obvious notation
  $H_2=\frac{1}{2}\pscal{\nu}{x^2+\xi^2}$.  One can decompose $H_2$
  into
  \begin{equation}
    H_2=\frac{1}{2}\pscal{\ell_1}{x^2+\xi^2} + \cdots +
    \frac{1}{2}\pscal{\ell_k}{x^2+\xi^2}
    \label{equ:oscillateurs}
  \end{equation}
  where $\ell_j\in\RM.\ZM^n$, $(\ell_1,\dots,\ell_k)$ are independent,
  and each Hamiltonian $\frac{1}{2}\pscal{\ell_j}{x^2+\xi^2}$ has a
  periodic flow.

  To show this, consider the orthogonal complement of the resonance module
  \[
  \R^\perp=\{\mu\in\ZM^n,\qquad \pscal{\alpha}{\nu}=0 \impliq
  \pscal{\alpha}{\mu}=0 \quad \forall \alpha\in\ZM^n\}.
  \]
  Let $(u^1,\dots,u^k)$ be a $\ZM$-basis of $\R^\perp$.  One can view
  the $\QM$-module $\QM.\R=\R\otimes\QM$ as a $\QM$-vector space,
  endowed with the $\QM$-scalar product induced by
  $\pscal{\cdot}{\cdot}$.  Then
  $(\R\otimes\QM)^\perp=\R^\perp\otimes\QM$ and, by density or some
  algebraic argument, $(\R\otimes\RM)^\perp=\R^\perp\otimes\RM$.
  Therefore $\nu\in\R^\perp\otimes\RM$, so we have $k\geq 1$ and one
  can decompose
  \[
  \nu=\sum_{j=1}^k \lambda_ju^j, \qquad \lambda_j\in\RM .
  \]
  We define $\ell_j:=\lambda_ju^j$. Since $u^j$ has integer
  coefficients, it is clear that the Hamiltonian
  $\frac{1}{2}\pscal{u^j}{x^2+\xi^2}$ has a periodic flow.

  From lemma~\ref{lemm:commutation} if $\beta,\gamma$ are
  multi-indices in $\NM^n$ and $z_j=x_j+i\xi_j$, then
  \begin{equation}
    \ad{H_2}(z^\beta\bar{z}^\gamma)= 
    \pscal{\gamma-\beta}{\nu}z^\beta\bar{z}^\gamma .
    \label{equ:commutation}
  \end{equation}
  But if $\gamma-\beta\in\R$, then
  \[
  \gamma-\beta \in ((\R\otimes\QM)^\perp)^\perp =
  (\R^\perp\otimes\QM)^\perp,
  \]
  Hence
  \[
  \forall j,\qquad \pscal{u_j}{\gamma-\beta}=0 .
  \]
  In other words $z^\beta\bar{z}^\gamma$ commutes with each term in
  the decomposition~\eqref{equ:oscillateurs}. Therefore any polynomial
  in $\ker\ad{H_2}$ commutes with all
  $\frac{1}{2}\pscal{\ell_j}{x^2+\xi^2}$'s, and thus is invariant
  under the $\T^k$ action they generate.
\end{demo}

\begin{rema}
  Given a point in $\RM^{2n}$, its orbit under the flow of $H_2$ is
  contained in the $\T^k$ orbit of that point. Actually a small
  variant of the proof shows that the inverses of the primitive
  periods of each periodic Hamiltonian
  $\frac{1}{2}\pscal{\ell_j}{x^2+\xi^2}$ possess no resonance
  relation, and hence the $H_2$-orbit is in fact \emph{dense} in the
  $\T^k$-orbit.
\end{rema}
\begin{coro}
  \label{coro:resonances}
  \begin{enumerate}
  \item If there is no resonance relation (\emph{ie.} $k=n$) then any
    element of $\E$ commuting with $H_2$ is of the form
    \[
    K = f(x_1^2+\xi_1^2,\dots,x_n^2+\xi_n^2;\h),
    \]
    for a formal series $f$ in $n+1$ variables.
  \item More generally if we let
    \[
    r=\inf\{\abs{\alpha}; \quad \alpha\in\ZM^n, \alpha\neq 0,
    \pscal{\alpha}{\nu}=0\}\quad\in\NM^*\cup\{\infty\}
    \]
    then any element of $\E$ commuting with $H_2$ is of the form
    \[
    K = f_r(x_1^2+\xi_1^2,\dots,x_n^2+\xi_n^2;\h) + R_r,
    \]
    where $R_r\in\O_r$ and $f_r(u;\h)$ is a polynomial in $(u;\h)$ of
    degree at most $[(r-1)/2]$.
  \end{enumerate}
\end{coro}
Of course this corollary follows even more obviously from
lemma~\ref{lemm:commutation} alone, since the monomials
$z^\alpha\bar{z}^\beta$ with $\abs{\alpha}+\abs{\beta}<r$ will commute
with $H_2$ only if $\alpha=\beta$. Hence they admit the form
$z^\alpha\bar{z}^\alpha=\prod(x_i^2+\xi_i^2)^{\alpha_i}$.

\section{The semiclassical Birkhoff normal form}
\label{sec:semi}
The goal of this section is to show how the formal Birkhoff normal
form can be transformed into a more usable semiclassical statement
involving spectral estimates. To make the proof more transparent, it
is enlightening to separate some statements which are independent of the
normal form construction, and which we believe have their own
interest.

In all the article we use the following notation. If $P$ is a
self-adjoint operator on some Hilbert space, $P$ bounded from below,
then the increasing sequence of eigenvalues below the essential
spectrum is denoted by
$\lambda^P_1\leq\lambda^P_2\leq\cdots\leq\lambda^P_j\leq\cdots$.  If
$I$ is a borelian of $\RM$, the spectral projector of $P$ on $I$ is
denoted by $\Pi_I^P$. If $P$ is a semiclassical pseudo-differential
operator, then of course $\lambda^P_j=\lambda^P_j(\h)$ and
$\Pi^P_I=\Pi^P_I(\h)$ also depend on $\h$.
 
\subsection{Semi-excited spectrum and Taylor expansions}

Let $X$ be either a compact manifold of dimension $n$ equipped with a
smooth density or $X=\RM^n$. We will deal with semiclassical
pseudo-differential operators on $X$ in the usual sense, as follows.
Let $d$ and $m$ be real numbers. When $X=\RM^n$, let $S^d(m)=S^d(m,X)$
the set of all families $(a(\cdot;\h))_{\h\in (0,1]}$ of functions in
$\Cinf(T^*X)$ such that
\begin{equation}
  \forall \alpha\in\NM^n, \quad \abs{\partial^\alpha_{(x,\xi)}
    a(x,\xi;\h)} \leq C_\alpha\h^d(1+\abs{x}^2+\abs{\xi}^2)^{\frac{m}2},
  \label{equ:symboles}
\end{equation}
for some constant $C_\alpha>0$, uniformly in $\h$.  Then $\Psi^d(m,X)$
is the set of all (unbounded) linear operators $A$ on $L^2(X)$ that
are $\h$-Weyl quantisations of symbols $a\in S^d(m)$~:
\[
(Au)(x) = (Op^w_\h(a)u)(x) =
\frac{1}{(2\pi\h)^n}\int_{\RM^{2n}}\!\!\!e^{\frac{i}{\h}\pscal{x-y}\xi}a({\textstyle\frac{x+y}2},\xi;\h)
u(y)\abs{dyd\xi}.
\]
The number $d$ in~(\ref{equ:symboles}) is called the order of the
operator. Unless specified, it will always be zero here.  In case $X$
is a compact manifold with a smooth density, $\Psi^d(m,X)$ is the set
of operators on $L^2(X)$ that are a locally finite sum $P=\sum_\beta
P_\beta+R$, where for each $\beta$ there is a open set $U_\beta\subset
X$ equipped with a chart $U_\beta\fleche \RM^n$ through which
$P_\beta\in\Psi^d(m,\RM^n)$, and $R$ is an integral operator whose
Schwartz kernel is $\O(\h^\infty)$ in the $\Cinf$ topology. Thus, if
$X$ is a compact riemannian manifold, $\Delta$ the corresponding
Laplacian, and $V\in\Cinf(X)$, the Schr{\"o}dinger operator
$P=-\frac{\h^2}2\Delta + V$ is a good candidate, of order zero. In
case $X=\RM^n$, the Schr{\"o}dinger operator is admissible whenever $V$
has at most a polynomial growth.

Let us denote $\Psi(m,X)=\cup_{d\in\ZM}\Psi^d(m,X)$,
$\Psi^d(X)=\cup_{m\in\ZM}\Psi^d(m,X)$, and
$\Psi(X)=\cup_{m\in\ZM}\Psi(m,X)$ .  We shall use in this article the
standard properties of such pseudo-differential operators.  In
particular the composition sends $\Psi(m,X)\times\Psi(m',X)$ to
$\Psi(m+m',X)$. Moreover all $P\in\Psi(0,X)$ are bounded:
$L^2(X)\fleche L^2(X)$, uniformly for $0<\h\leq 1$.

If $P$ has a real-valued Weyl symbol, then it is a symmetric operator
on $L^2$ with domain $\Cinf_0(X)$. If its principal symbol is bounded
from below then we use the Friedrichs self-adjoint extension, and we
will identify $P$ with this extension. Actually, if the Weyl symbol is
real and $p$ is elliptic at infinity (\emph{i.e.} if there exists
$m\in\RM$ and a constant $C>0$ such that $P\in \Psi(m,X)$ and its
principal symbol $p$ satisfies $\abs{p(x,\xi)}\geq
\frac1C(\norm{x}^2+\norm{\xi}^2)^{m/2}$ for
$\norm{x}^2+\norm{\xi}^2\geq C$), then $P$ is essentially self-adjoint
(see for instance~\cite[proposition 8.5]{dimassi-sjostrand}). But we
won't use this result here.

Finally, when $P\in\Psi(m,X)$ is self-adjoint and $f\in\Cinf_0(\RM)$,
then $f(P)\in\cap_{m'}(\Psi(m',X))$. See for
instance~\cite{dimassi-sjostrand},~\cite{robert},
or~\cite{colin-livre} for details.  In this work all
pseudo-differential operators are assumed to admit a classical
asymptotic expansion in integer powers of $\h$.

\begin{theo}
  \label{theo:taylor}
  Let $P$ and $Q$ be two semiclassical pseudo-differential operators
  in $\Psi^0(X)$ such that
  \begin{itemize}
  \item at $z_0\in T^*X$, the principal symbols $p$ and $q$ take their
    minimal value $p(z_0)=q(z_0)=0$, this minimum is reached only at
    $z_0$ and is non-degenerate;
  \item there exists $E_\infty>0$ such that $\{p\leq E_\infty\}$ and
    $\{q\leq E_\infty\}$ are compact.
  \end{itemize}
  Suppose that, in some local coordinates near $z_0$, the total
  symbols of $P$ and $Q$ admit the same Taylor expansion at $z_0$.
  Then there exists $E_0>0$, $\h_0>0$ and for each $N$ a constant
  $C_N>0$ such that for all $(\h,E)\in[0,\h_0]\times[0,E_0]$
  \[
  \lambda^P_j\leq E \textup{ or } \lambda_j^Q\leq E \impliq
  \abs{\lambda_j^P-\lambda_j^Q}\leq C_N(E^N+\h^N) .
  \]
\end{theo}
Before entering the proof of the theorem, we recall an elementary
consequence of the minimax theorem.
\begin{lemm}
  \label{lemm:comparaison}
  Let $A$ and $B$ be two self-adjoint operators on $\H$, both bounded
  from below. Suppose there exists an interval $I=(-\infty, E]$ and
  $C>0$ such that $\Pi_I^B(\H)\subset\textup{Dom}(A)$ and
  \[
  \norm{(A-B)\Pi^B_I}\leq C .
  \]
  Then for all $j$ such that $\lambda^B_j\leq E$ one has
  \[
  \lambda^A_j\leq \lambda^B_j + C .
  \]
\end{lemm}
\begin{demo}
  Let $\lambda^B_j\leq E$ and
  $\F^B_j=\Pi^B_{(-\infty,\lambda^B_j]}(\H)$ the eigenspace associated
  to the eigenvalues below or equal to $\lambda^B_j$. Let
  $\phi\in\F^B_j$, of norm $1$; by hypothesis one has
  \[
  \norm{A\phi-B\phi}\leq C .
  \]
  Hence by Cauchy-Schwarz
  $\abs{\pscal{A\phi}{\phi}-\pscal{B\phi}{\phi}}\leq C$. Therefore
  \[
  \pscal{A\phi}{\phi}\leq \pscal{B\phi}{\phi} + C .
  \]
  Since
  $\lambda^B_j=\sup_{\phi\in\F^B_j,\norm{\phi}=1}\pscal{B\phi}{\phi}$,
  one gets
  \[
  \sup_{\phi\in\F^B_j,\norm{\phi}=1}\pscal{A\phi}{\phi}\leq
  \lambda^B_j+C .
  \]
  Now since $\F^B_j$ has dimension $j$, the minimax formula
  \[
  \lambda^A_j = \inf_{\F\subset\textup{Dom}(A), \dim
    \F=j}\left(\sup_{\phi\in\F,\norm{\phi}=1}
    \pscal{A\phi}{\phi}\right)
  \]
  implies that $\lambda^A_j\leq\lambda^B_j+C$.
\end{demo}

When dealing with manifolds $X\neq\RM^n$, we shall need a refinement
of the lemma, as follows.
\begin{lemm}
  \label{lemm:comparaison2}
  Let $A$ and $B$ be two self-adjoint operators acting respectively on
  the Hilbert spaces $\H'$ and $\H$, both bounded from below. Suppose
  there exists a bounded operator $U:\H\fleche\H'$, an interval
  $I=(-\infty, E]$ and constants $C>0$, $c\in(0,1)$ such that
  $U\Pi_I^B(\H)\subset\textup{Dom}(A)$ and
  \[
  \norm{(U^*AU-B)\Pi^B_I}\leq C
  \]
  and
  \[
  \norm{U^*U\Pi_I^B-\Pi^B_I}\leq c
  \]
  Then for all $j$ such that $\lambda^B_j\leq E$ one has
  \[
  \lambda^A_j\leq (\lambda^B_j + C)(1+\frac{c}{1-c}) .
  \]
\end{lemm}
\begin{demo}
  Using the same notation as in the proof of
  lemma~\ref{lemm:comparaison}, we deduce from the first hypothesis
  that
  \[
  \pscal{A U \phi}{U \phi}\leq \pscal{B\phi}{\phi} + C ,
  \]
  while the second yields
  \[
  \abs{\norm{ U \phi }^2 - 1} \leq c
  \]
  Hence $\norm{U\phi}^2\geq 1-c>0$. Therefore
  \[
  \frac{ \pscal{A U \phi}{U \phi}}{ \norm{U \phi} ^2} \leq
  (\pscal{B\phi}{\phi} + C)( 1+ \frac{c}{1-c} ).
  \]
  Moreover $U:\F_j^B\fleche\H'$ is injective and hence
  $\dim(U\F_j^B)=j$. We conclude as in the proof of
  lemma~\ref{lemm:comparaison}.
\end{demo}

Finally, for the proof of theorem~\ref{theo:taylor} it will be very
convenient to use a generalisation of a well-known microlocalisation
result of~\cite{sjostrand-semi} for which, using the above lemmas, we
give a new and simple proof. Recall that we say that a
pseudo-differential operator $P\in\Psi(X)$ microlocally vanishes at a
point $z\in T^* X$ when in some local coordinates its full Weyl symbol
vanishes at $z$.
\begin{lemm}
  \label{lemm:microlocal}
  Let $P\in\Psi(X')$ and $Q\in\Psi(X)$ be self-adjoint semiclassical
  pseudo-differential operators, with principal symbols $p$ and $q$.
  Assume there exists a bounded operator $U:L^2(X)\fleche L^2(X')$,
  compact subsets $D\subset T^*X$, $D'\subset T^*X'$, and an interval
  $I=(-\infty,E]$ (with $D$, $D'$ and $E$ being independent of $\h$)
  such that
  \begin{enumerate}
  \item $p^{-1}(I)$ (respectively $q^{-1}(I)$) is contained in the
    interior of $D$ (respectively $D'$);
  \item $U^*PU-Q$ and $U^*U-\Id$ are pseudo-differential operators that
    microlocally vanish in $D$;
  \item $P-UQU^*$ and $UU^*-\Id$ are pseudo-differential operators
    that microlocally vanish in $D'$;
  \end{enumerate}
  Then there exists $\h_0>0$ and a positive sequence $(C_N)_{N\in\NM}$
  (depending on $E$) such that, for all $j$ and $\h<\h_0$ such that
  $\lambda_j^Q\leq E$ (or $\lambda_j^P\leq E$), one has:
  \[
  \abs{\lambda_j^Q-\lambda_j^P}\leq C_N\h^N .
  \]
\end{lemm}
\begin{rema}
  In most situations $U$ will be a Fourier integral operator. When
  $X=X'$ and $P,Q\in\Psi(0,X)$, the result was proved
  in~\cite[proposition 2.2]{sjostrand-semi}, using the Kato distance
  to handle the spectral perturbation. Our proof here, using the
  minimax, looks simpler, but the idea is essentially the same. A
  small additional argument is needed to handle $X\neq X'$. We give
  the full proof here for the sake of completeness.
\end{rema}
\begin{demo}
  First let $E'>E$ such that the hypothesis $(1.)$ still holds when
  $E$ is replaced by $E'$.

  The hypothesis $(1.)$ ensures that the spectra of $P$ and $Q$ are
  strictly bounded from below by some positive constant $E_0$,
  independent of $\h$.  Moreover it is well known that it also implies
  that the intersections with $I'=(-\infty,E']$ of these spectra are
  discrete.  One can check this as follows.
  
  Let $E_1>E'$ be independent of $\h$ and such that
  $q^{-1}([E_0,E_1])$ is contained in the interior of $D$.  Let
  $f\in\Cinf_0(\RM)$ such that
  \begin{enumerate}
  \item $f=1$ on $[\lambda_1^Q,E']$ ;
  \item $f=0$ outside of $[E_0,E_1]$ .
  \end{enumerate}
  By pseudo-differential functional calculus (see for
  instance~\cite[th{\'e}or{\`e}me III-11]{robert}), $f(Q)$ is a
  pseudo-differential operator which belongs to the trace class and
  hence is compact. This entails that the spectral projector of $Q$
  onto $[E_0,E']$ is compact as well, thus proving the discreteness of
  $\Sp(Q)$ in $[E_0,E']$. Of course the same argument applies to $P$.
  See also~\cite{helffer-robert-calcul} for more details.

  What's more, the functional calculus also ensures that $f(Q)$
  microlocally vanishes outside $q^{-1}([E_0, E_1])$.  By symbolic
  calculus one has
  \[
  \norm{(Q-U^*PU)f(Q)}=\ohb .
  \]
  Therefore $\textup{Dom}(Qf(Q))=\textup{Dom}(U^*PUf(Q))$, in the
  sense of Friedrichs extensions. But for all
  $u\in\Pi^Q_{I'}(\H)$, $u=f(Q)u$. Hence $Uu\in\textup{Dom}(P)$ and
  \[
  \norm{(Q-U^*PU)u} = \norm{(Q-U^*PU)f(Q)u} = \ohb\norm{u}.
  \]
  This shows that there exists a positive sequence $(C_N)_{N\in\NM}$
  such that for all $N$, $\norm{(Q-U^*PU)\Pi^Q_{I'}}\leq C_N\h^N$.
  
  Similarly, $\norm{(U^*U-\Id)f(Q)}=\ohb$ and hence
  \[
  \norm{(U^*U-\Id)\Pi^Q_{I'}}\leq c_N\h^N,
  \]
  for a positive sequence $(c_N)_{N\in\NM}$.

  Applying now lemma~\ref{lemm:comparaison2} we get, for all $j$ such
  that $\lambda_j^Q\leq E'$, the inequality
  $\lambda_j^P\leq(\lambda_j^Q+C_N\h^N)(1+\frac{c_N\h^N}{1-c_N\h^N})$.
  In particular for $\h$ small enough we always reach $\lambda_j^P\leq
  E'$ whenever $\lambda^Q_j\leq E$.

  Interchanging the roles of $P$ and $Q$ we obtain as well
  $\norm{(P-UQU^*)\Pi^P_{I'}}\leq C_N\h^N$ and
  $\norm{(UU^*-\Id)\Pi^P_{I'}}\leq c_N\h^N$ (with perhaps a
  modification of $C_N$ and $c_N$). Hence a new application of
  lemma~\ref{lemm:comparaison2} yields
  $\lambda^Q_j\leq(\lambda^P_j+C_N\h^N)(1+\frac{c_N\h^N}{1-c_N\h^N})$,
  uniformly for all $j$ such that $\lambda^P_j\leq E'$. This shows
  that
  \[
  \abs{\lambda_j^Q-\lambda_j^P}\leq (C'_N+E')\h^N .
  \]
  as soon as $\lambda_j^Q\leq E$. Swapping again the roles of $P$ and
  $Q$ we obtain the final result.
\end{demo}

\begin{demo}[of theorem~\ref{theo:taylor}]
  The result of the theorem will be denoted as the property
  $\P(P,Q,\h_0,E_0,C_N)$. It is easy to see that if
  $\P(P,Q,\h_0,E_0,C_N)$ and $\P(Q,R,\h_0',E_0',C_N')$ hold, then
  $\P(P,R,\h_0'',E_0'',C_N'')$ will hold with suitably chosen
  constants $\h_0'',E_0'',C_N''$.

  \paragraph{1. ---} We use this transitivity property to
  microlocalise the problem in a compact subset of $T^*X$. Let $\Phi$
  be a pseudo-differential operator that is microlocally equal to the
  identity on a neighbourhood of the compact
  \[
  D:=\{p\leq E_\infty\}\cup\{q\leq E_\infty\}
  \]
  and microlocally vanishes outside a compact of $T^*X$. We may assume
  also that its principal symbol $\phy$ satisfies $0\leq \phy\leq 1$.
  Then consider the operator $P'=\Phi P+2E_\infty(\textup{Id}-\Phi)$.
  Its principal symbol $p'$ satisfies $p'=p$ for $p\leq E_\infty$ and
  $p'>E_\infty$ as soon as $p>E_\infty$. Hence $P'$ satisfies the
  hypothesis of the theorem as well. Using that $P'$ is microlocally
  equal to $P$ on $D$, we apply lemma~\ref{lemm:microlocal} to $P$ and
  $P'$ with $U=\Id$ and an energy $E=E_0<E_\infty$ such that $\{p\leq
  E\}$ is contained in the interior of $D$. Then, since $\h^N\leq
  E^N+\h^N$, we obtain $\P(P,P',\h_0,E_0,C_N)$. Using the same trick
  for $Q$, we construct $Q'$ with $\P(Q,Q',\h_0',E_0',C_N')$ for new
  constants $h_0',E_0',C_N'$.

  Thus by transitivity we are reduced to prove the theorem for $P'$
  and $Q'$.

  \paragraph{2.---} We compare now $P'$ and $Q'$. Notice that
  $R:=P'-Q'=\Phi(P-Q)$ microlocally vanishes outside a compact subset
  of $T^*X$. By hypothesis the Weyl symbols of $P'$ and $Q'$ near
  $z_0$ have the same Taylor expansion. Hence the symbol of $R$ is
  flat at $z_0$. By symbolic calculus we can construct a
  pseudo-differential operator $S_N$ such that
  \begin{equation}
    R=S_N (P')^N + \ohb ,
    \label{equ:division}
  \end{equation}
  and $S_N$, as $R$ does, microlocally vanishes outside a compact of
  $T^*X$. This implies that $S_N$ is bounded for $\h\leq 1$ by a
  constant independent of $\h$.

  Hence~\eqref{equ:division} implies, for all $E>0$ and $\h\leq 1$,
  the following estimate
  \begin{equation}
    \norm{R\Pi^{P'}_{[-E,E]}}\leq D_N(E^N+\h^N) .
    \label{equ:majore}
  \end{equation}
  We claim that there is a positive sequence $(C_N)$ such that
  \[
  \norm{R\Pi^{P'}_{(-\infty,E]}}\leq C_N(E^N+\h^N) .
  \]
  Indeed let $-E_\textup{min}$ be the bottom of the spectrum of $P'$.
  If $E_\textup{min}\leq E$ then $\Pi^{P'}_{(-\infty,
    E]}=\Pi^{P'}_{[-E,E]}$ and the formula follows
  from~\eqref{equ:majore} with $C_N=D_N$. If $E_\textup{min}>E$ then
  \[
  \norm{R\Pi^{P'}_{(-\infty,E]}} \leq
  \norm{R\Pi^{P'}_{[-E_\textup{min},E_\textup{min}]}} \leq
  D_N(E_\textup{min}^N+\h^N) .
  \]
  But by Garding's inequality there is a constant $C>0$ such that
  $E_\textup{min}\leq C\h$. Hence
  \[
  \norm{R\Pi^{P'}_{(-\infty,E]}}\leq D_N(C^N+1)\h^N\leq C_N(E^N+\h^N)
  \]
  with $C_N=D_N(C^N+1)$. Thus the claim is proved.

  Now let $E_0>0$ be such that the spectrum of $P'$ is discrete in
  $(-\infty, E_0]$. Lemma~\ref{lemm:comparaison} ensures us that if
  $(\h,E)\in[0,1]\times[0,E_0]$ then
  \[
  \lambda^{P'}_j\leq E \impliq
  \lambda^{Q'}_j\leq\lambda^{P'}_j+C_N(E^N+\h^N) .
  \]

  Finally, as in the proof of lemma~\ref{lemm:microlocal}, we may
  interchange the roles of $P'$ and $Q'$ to obtain
  $\P(P',Q',\h_0',E_0',C_N')$, for some new positive constants
  $\h_0'$, $E_0'$ and $C_N'$.
\end{demo}

\subsection{Using the formal Birkhoff normal form}
In this section we consider a pseudo-differential operator
$P\in\Psi(X)$ fulfilling the hypothesis of theorem~\ref{theo:taylor},
transform it into a pseudo-differential operator on $\RM^n$, take its
Taylor series at $z_0$, apply the formal Birkhoff normal form, and
finally construct a new pseudo-differential operator $Q$ commuting
with a harmonic oscillator $H_2$. $Q$ is compared with $P$ using
lemma~\ref{lemm:microlocal} and theorem~\ref{theo:taylor}, hence
reducing the spectral study of $P$ to that of an effective Hamiltonian
on some eigenspace of the harmonic oscillator.

\ouf

Thus, first of all, we transfer the spectral problem to $\RM^n$.
\begin{lemm}
  \label{lemm:Rn}
  Let $P\in\Psi(X)$ satisfy the hypothesis of
  theorem~\ref{theo:taylor} at a point $z_0\in T^*X$. Then there
  exists a pseudo-differential operator $Q\in\Psi(\RM^n)$ satisfying
  the hypothesis of theorem~\ref{theo:taylor} at the origin
  $0\in\RM^{2n}$, some constants $E>0$ and $\h_0>0$, and a positive
  sequence $(C_N)_{N\in\NM}$ such that, for all $j$ and $\h<\h_0$ such
  that $\lambda_j^Q\leq E$ (or $\lambda_j^P\leq E$), one has:
  \[
  \abs{\lambda_j^Q-\lambda_j^P}\leq C_N\h^N .
  \]
  Moreover there exist local coordinates near $z_0$ in which the full
  Weyl symbol of $P$ is exactly the Weyl symbol of $Q$.

\end{lemm}
\begin{demo}
  Let $(x,\xi)$ be canonical coordinates on a neighbourhood $\Omega$
  of $z_0=(x_0,\xi_0)$ coming from local coordinates $x$ on $X$ and
  let $U=U(\h)$ be the integral operator with Schwartz kernel
  \[
  U(x,y) = \frac{1}{(2\pi\h)^n}\int
  e^{\frac{i}{\h}\pscal{x-y}{\xi}}\phy(y,x,\xi)d\xi,
  \]
  where $x\in X$, $y\in\RM^n$ and $\phy\in\Cinf_0(\RM^n\times\Omega)$
  with $\phy(y,x,\xi)\equiv 1$ in a neighbourhood of $(y_0,z_0)$. Then
  $U:L^2(\RM^n)\fleche L^2(X)$ is bounded and $U^*PU\in\Psi(\RM^n)$
  and has, when expressed in the coordinates $(x,\xi)$, the same Weyl
  symbol as $P$. This follows from the fact that, restricted to
  $\Omega$ and expressed in these coordinates, $U$ is simply a
  compactly supported pseudo-differential operator microlocally equal
  to the identity near the origin. (From a more geometrical viewpoint
  $U$ is actually a Fourier integral operator associated to the
  symplectomorphism $T^*X\fleche\RM^{2n}$ defined by the canonical
  coordinates $(x,\xi)$).  Since the principal symbol of $U^*PU$ has a
  local non-degenerate minimum at the origin, an easy
  pseudo-differential partition of unity will modify $U^*PU$ outside a
  microlocal neighbourhood $D$ of the origin in such a way that its
  principal symbol will satisfy the global hypothesis of
  theorem~\ref{theo:taylor}. Let $Q$ be the modified operator. Then
  for $E>0$ small enough, the hypothesis of
  lemma~\ref{lemm:microlocal} are fulfilled (with $D'$ identified with
  $D$ thanks to the local coordinates in $\Omega$).  This lemma gives
  exactly the desired spectral result.
\end{demo}
\begin{rema}
  We shall not use the fact that $P$ and $Q$ have the same Weyl
  symbols in some coordinates. We only retain the geometrical fact
  that their principal and sub-principal symbols are symplectomorphic.
\end{rema}

Using the canonical coordinates of $T^*\RM^n$ we introduce the space
$\E$ as in section~\ref{sec:formal}. Using a Borel resummation, one
can always quantise an element in $L\in\E$ into a pseudo-differential
operator in $\Psi(\RM^n)$ whose Weyl symbol has a Taylor series giving
back the initial series in $\E$. Moreover we can arbitrarily extend
the pseudo-differential operator to vanish microlocally far from
$z_0$. With a slight abuse with respect to the standard notation, we
shall in this section denote by $Op^W(L)$ such a pseudo-differential
operator; and for any $Q\in\Psi(\RM^n)$ we denote by
$\sigma_W(Q)\in\E$ the Taylor series at $z_0$ of the Weyl symbol of
$Q$.
 
Let $Q\in\Psi(\RM^n)$ satisfy the hypothesis of
theorem~\ref{theo:taylor} at $z_0=0\in\RM^{2n}$.  Consider the Taylor
series $[Q]=\sigma_W(Q)$. Since $z_0$ is a non-degenerate minimum for
$p$, one has
\[ [Q]=\h H_0 + H_2 + L
\]
where $L\in\O(3)$, $H_2\in\D_2$ is elliptic, and $H_0\in\RM$ is the
value at the origin of the sub-principal symbol of $Q$. Applying the
formal Birkhoff normal form of theorem~\ref{theo:formal} we obtain
$[A]$ and $[K]$ in $\O_3$ such that $e^{i\h^{-1}\ad{[A]}}(H_2+[L]) =
H_2+[K]$.  Consider the operators $A=Op^W([A])$ and $K=Op^W([K])$ (so
that $\sigma_W(A)=[A]$ and $\sigma_W(K)=[K])$. Now $e^{i\h^{-1}A}$ is
a Fourier integral operator and by Egorov's theorem
$e^{i\h^{-1}A}Qe^{-i\h^{-1}A}\in\Psi(\RM^n)$.
\begin{lemm}
  The Taylor expansion at $z_0$ of the Weyl symbol of
  $e^{i\h^{-1}A}Qe^{-i\h^{-1}A}$ is $\h H_0 + H_2+[K]$.
\end{lemm}
\begin{demo}
  Since $A$ is bounded
  $e^{i\h^{-1}A}Qe^{-i\h^{-1}A}=\exp(\ad{i\h^{-1}A})Q$. Expanding the
  exponential in the right-hand side using Taylor's formula with
  integral remainder, one gets
  \[
  e^{i\h^{-1}A}Qe^{-i\h^{-1}A} = \sum_{j=0}^N
  \frac{1}{j!}(\ad{i\h^{-1}A})^jQ +
  \frac{1}{N!}(\ad{i\h^{-1}A})^{N+1}R_N ,
  \]
  where
  \[
  R_N=\int_0^1(1-t)^N(\exp(t\ad{i\h^{-1}A})Q)dt .
  \]
  By definition of the Lie algebra structure of $\E$,
  $\sigma_W(\sum_{j=0}^N \frac{1}{j!}(\ad{i\h^{-1}A})^jQ)$ is
  precisely $\sum_{j=0}^N \frac{1}{j!}(\ad{i\h^{-1}[A]})^j[Q]=\h
  H_0+H_2+[K]+\O_{N+1}$. Thus we need to prove that
  $\sigma_W(\frac{1}{N!}(\ad{i\h^{-1}A})^{N+1}R_N)\in\O_{N+1}$, and
  for this it suffices to show that the Weyl symbol of $R_N$ is
  bounded near $z_0$, uniformly in $\h$. Indeed, its Taylor series
  would then be in $\O_0$ and we would conclude using the fact that
  $\ad{i\h^{-1}[A]}$ sends $\O_j$ to $\O_{j+1}$.  But by Egorov's
  theorem, $\exp(t\ad{i\h^{-1}A})Q$ is a pseudo-differential operator
  of order zero, uniformly in $t\in[0,1]$. Integrating over $t$ we get
  that $R_N$ is indeed of order zero.
\end{demo}

\begin{prop}
  \label{prop:weyl}
  For any compact $D\subset\RM^{2n}$ containing the origin in its
  interior, there exists a pseudo-differential operator
  $K\in\Psi(\RM^n)$, microlocally vanishing outside $D$, such that
  \begin{itemize}
  \item $[\hat{H}_2,K]=0$;
  \item The Weyl symbols of $e^{i\h^{-1}A}Qe^{-i\h^{-1}A}$ and $\h
    H_0+\hat{H}_2+K$ have the same Taylor expansion at the origin;
  \item $\h H_0+\hat{H}_2+K$ satisfies the hypothesis of
    theorem~\ref{theo:taylor} (with $X=\RM^n$ and $z_0=0$).
  \end{itemize}
\end{prop}

\begin{demo}
  We use lemma~\ref{lemm:torus}, \emph{ie.} the fact that there exists
  a $\T^k$ action on $\RM^{2n}_{\{x,\xi\}}$ (with $k\geq 1$) such that
  \begin{equation}
    \ker \ad{H_2} = \E^{\T^k} .
    \label{equ:torus}
  \end{equation}
  Let $\tilde{K}$ be a compactly supported Borel resummation of $[K]$.
  Let $\bar{K}$ be the average of $\tilde{K}$ under the $\T^k$-action.
  Then, since $[K]\in\E^{\T^k}$, $\sigma_W(\bar{K})=[K]$. Hence
  $\sigma_W(\bar{K})\in\E^{\T^k}$ and by Weyl quantisation,
  $K:=Op^W(\bar{K})$ commutes with $\hat{H_2}$. (Recall that
  commutation is preserved by Weyl quantisation here because $H_2$ is
  quadratic.)

  Since $[K]\in\O_3$, the last point of the proposition is
  automatically satisfied if one chooses the support of the Borel
  resummation to be close enough to the origin.
\end{demo}

From the proposition, we deduce:
\begin{coro}
  \label{coro:fnq0}
  The operators $Q$ and $\h H_0+\hat{H}_2+K$ have equivalent spectra
  in the sense of theorem~\ref{theo:taylor}.
\end{coro}
\begin{demo}
  Apply theorem~\ref{theo:taylor} to $e^{i\h^{-1}A}Qe^{-i\h^{-1}A}$
  and $\h H_0+\hat{H}_2+K$.
\end{demo}

We are now in position to state and prove the quantum Birkhoff normal
form, which reduces the spectral analysis of $P$ in the semiclassical
regime to that of a reduced Hamiltonian, $K$, acting on some
eigenspace of $H_2$ of finite dimension (but growing as $\h$ decreases
or $E$ increases).

\begin{theo}
  \label{theo:FNB}
  Let $P\in\Psi(X)$ be a semiclassical self-adjoint
  pseudo-differential operator of order zero such that
  \begin{itemize}
  \item at $z_0\in T^*X$, the principal symbol $p$ takes its minimal
    value $p(z_0)=0$, this minimum is reached only at $z_0$ and is
    non-degenerate;
  \item there exists $E_\infty>0$ such that $\{p\leq E_\infty\}$ is
    compact.
  \end{itemize}
  Let $H_0$ be the value at $z_0$ of the sub-principal symbol of $P$.
  Then there exists a harmonic oscillator $H_2$ on $\RM^{2n}$ (ie. an
  elliptic element of $\D_2$ in the terminology of
  section~\ref{sec:formal}), and for any compact domain
  $D\subset\RM^{2n}$ containing the origin in its interior there
  exists a pseudo-differential operator $K\in\Psi(\RM^n)$ of order
  zero such that
  \begin{itemize}
  \item $[K,\hat{H}_2]=0$;
  \item $K$ vanishes microlocally outside of $D$;
  \item $\sigma_W(K)\in\O_3$,
  \end{itemize}
  and for each $\eta>0$ there exists $E_0>0$, $\h_0>0$ and for each
  $N$ a constant $C_N>0$ such that for all
  $(\h,E)\in[0,\h_0]\times[0,E_0]$,
  \[
  \left(\lambda^P_j\leq E \textup{ or } \lambda_j^{Q}\leq E - \h
    H_0\right) \impliq \abs{\lambda_j^P-\lambda_j^{Q} - \h H_0}\leq
  C_N(E^N+\h^N) ,
  \]
  where
  \[
  Q=Q((1+\eta)E) := (\hat{H_2}+K)_{\restr
    \Pi^{\hat{H}_2}_{(-\infty,(1+\eta)E]}(L^2(\RM^n))} .
  \]
\end{theo}
\begin{demo}
  We first apply lemma~\ref{lemm:Rn} which allows us to assume that
  $P\in\Psi(\RM^n)$ and $z_0=0$.

  In view of corollary~\ref{coro:fnq0} the theorem already holds with
  $Q((1+\eta)E)$ replaced by $\hat{H_2}+K$. But for $E$ and $D$ small
  enough the spectra of $Q((1+\eta)E)$ and $\hat{H_2}+K$ are exactly
  the same. To see this, we use the following estimate, which is due
  to the fact that $[K]\in\O_3$:
  \begin{equation}
    \exists \h_0>0, \exists C>0,\quad \forall E'>0,\forall \h<\h_0 \qquad
    \norm{K\Pi^{\hat{H_2}}_{(-\infty,E']}}\leq C{E'}^{\frac32} .
    \label{equ:estimeeK}  
  \end{equation}
  This estimate is a particular case of a more general result proven
  below (lemma~\ref{lemm:normes}).

  Let $E_c=4/9C^2$ (this is where the function $E'\mapsto
  E'-C{E'}^\frac32$ reaches its maximal value). Then without
  modifying~\eqref{equ:estimeeK} one can assume that $D$ is included
  in $H_2^{-1}((-\infty,E_c])$: for this one can replace $K$ by
  $Kf(\hat{H_2})$ where $f\in\Cinf_0(\RM)$ takes values in $[0,1]$, is
  equal to $1$ near the origin, and is supported inside
  $(-\infty,E_c]$. Then the estimate~\eqref{equ:estimeeK} can be
  improved as follows:
  \begin{equation}
    \label{equ:estimeeK2}
    \forall E'>0, \qquad
    \norm{K\Pi^{\hat{H_2}}_{(-\infty,E']}} \leq r(E') ,
  \end{equation}
  where $r(E'):=\min(C{E'}^\frac32,C{E_c}^\frac32)$.  For any $E'>0$
  and any operator $Q$ we use the notation
  \[
  \H^Q_{E'}:=\Pi^Q_{(-\infty,E']}(L^2(X)) .
  \]
  Also let $Q(\infty):=\hat{H_2}+K$.  Using that $\H^{Q(\infty)}_{E}$
  is stable by $\hat{H}_2$ and $\hat{H}_2=Q(\infty)-K$ we see that
  \begin{equation}
    \H^{Q(\infty)}_{E}\subset \H^{\hat{H}_2}_{\tilde{E}} ,
    \label{equ:hilbert}
  \end{equation}
  with $\tilde{E}=E+\norm{K\Pi^{Q(\infty)}_{E}}$.
  But~\eqref{equ:hilbert} implies $\norm{K\Pi^{Q(\infty)}_{E}}\leq
  r(\tilde{E})$. Inverting the function
  $\tilde{E}\mapsto\tilde{E}-r(\tilde{E})$ which by construction is
  strictly increasing for $\tilde{E}\geq 0$ we see
  that~\eqref{equ:hilbert} holds as soon as $\tilde{E}\geq
  E+C{E}^\frac32 +\O({E}^2)$, which is satisfied if
  $\tilde{E}=(1+\eta)E$, provided $E$ is small enough.
  Then~\eqref{equ:hilbert} says that the eigenvalues of $Q(\infty)$
  less than $E$ are the same as the eigenvalues less that $E$ of the
  restriction of $Q(\infty)$ to $\H^{\hat{H}_2}_{\tilde{E}}$. This
  restriction is precisely $Q((1+\eta)E)$.
\end{demo}

\section{The joint spectrum}
\label{sec:joint}

In the previous section, we have split $P$ into a harmonic oscillator
$\hat{H}_2$ and a commuting perturbation $K$ whose microsupport could
be arbitrarily small, provided one is only interested in sufficiently
small eigenvalues of $P$. As in~\cite{sjostrand-semi}, this can be
exploited to describe the so-called \emph{semi-excited states}, whose
energies are of order $\O(h^\gamma)$, $\gamma\in(0,1)$. Here we push
the analysis one step further, by introducing a second semiclassical
parameter associated to the high energy of the harmonic oscillator
$\hat{H}_2$ and which turns out to govern the study of the
perturbation $K$. In some sense the idea is to have a semiclassical
regime associated to the \emph{reduction} of $K$ by the $\hat{H}_2$
action. This view point is made explicit in the next section where we
shall assume that this action is periodic. For the moment the game is
to control $K$ in a semiclassical regime with two semiclassical
parameters...

So, let

\begin{equation}
  H_2(x,\xi)=\sum_{j=1}^n\frac{\nu_j}{2}(x_j^2+\xi_j^2),
  \label{equ:harmonic_osc}
\end{equation}
defined on $\RM^{2n}$, and let $\hat{H}_2$ be its Weyl quantisation,
acting on $L^2(\RM^n)$.

Let $K=K(\h)\in\Psi(\RM^n)$ a \pdo{} commuting with $\hat{H}_2$. Let
$\H_{E,\h}$ be the eigenspace of $\hat{H}_2$ for the eigenvalue $E$.
Our goal is to study the restriction of $K$ to $\H_{E,\h}$, in terms
of the parameters $E$ and $\h$, both in a neighbourhood of the origin.
More precisely, let $E_0>0$ and restrict the set of admissible $E$'s
to eigenvalues of $\hat{H}_2$ less that $E_0$, \emph{ie.} we consider
the set
\[
\{(E,\h)\in(0,E_0)\times(0,\infty) / \qquad
\exists\alpha\in\NM^n,\quad
E=\h(\textstyle\frac{\abs{\nu}}{2}+\pscal{\nu}{\alpha})\} .
\]

Instead of working with $(E,\h)$ we shall use the more convenient
scaling $(\epsilon,h)$ defined by
\begin{equation}
  \left\{
    \begin{array}{l}
      \epsilon(E,\h)=\sqrt{E}\\
      h(E,\h)=\h/E
    \end{array}\right. \qquad \epsilon\in(0,\epsilon_0], h^{-1}\in
  {\textstyle\frac{\abs{\nu}}{2}} + \pscal{\nu}{\NM^n} ,
  \label{equ:rescaling}
\end{equation}
with $\epsilon_0=\sqrt{E_0}$.
\begin{figure}[h]
  \centering
  \includegraphics[width=0.9\textwidth]{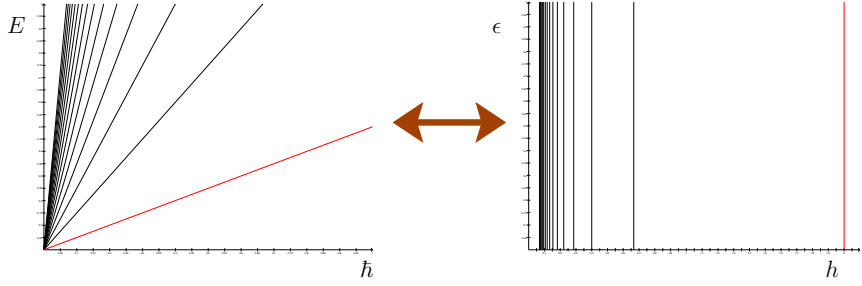}
  \caption{Allowed values for the parameters $(\h,E)$ or
    $(h,\epsilon)$.}
  \label{fig:parameters}
\end{figure}
Notice that the dimension of $\H_{E,\h}$ can be expressed as a
function of $h$ alone. Let $\H_h=\H_{1,h}$. We identify $\H_{E,\h}$
with $\H_h$ using the unitary transformation
\begin{equation}
  U_\epsilon : \H_{E,\h}\fleche \H_h, \qquad f(x)\mapsto
  \epsilon^{\frac{n}{2}}f(\epsilon x) ,
  \label{equ:scaling}
\end{equation}
The extension of $U_\epsilon$ to $L^2 ( \RM ^n)$ satisfies, for any
symbol $p$,
\begin{equation}
  U_\epsilon Op^w_\h(p(.,\h)) U_\epsilon^* = Op^w_{h}(q(.,h))
  ,
  \label{equ:scaling-weyl}
\end{equation}
where $q(x,\xi;h):=p(\epsilon x,\epsilon\xi;\epsilon ^2 h)$.

\begin{defi}
  \label{defi:Mlm}
  Let $K=K(\h)$ be a pseudo-differential operator in $\Psi(\RM^n)$.
  For any positive integers $\ell$, $m$ we shall say that
  $K\in\M_{\ell,m}$ if
  \begin{itemize}
  \item its Weyl symbol $k$ is of order $\O(\h^m)$ in a neighbourhood
    of $\{H_2\leq E_0\}$;
  \item the Taylor expansion of $k$ at the origin is in $\O_\ell$
    (with the notation of section~\ref{sec:formal}).
  \end{itemize}
\end{defi}
If $K(\h)$ commutes with $\hat{H}_2$, we denote by
$(K(\epsilon,h):\H_h\fleche\H_h)_{(\epsilon,h)}$ the family of
operators defined by
\[
K(\epsilon,h)=U_\epsilon K(\h)_{\restr \H_{E,\h}}U_\epsilon^* .
\]
The main technical result of this section is the following estimate
\begin{lemm}
  \label{lemm:normes}
  Let $K$ be a pseudo-differential operator commuting with
  $\hat{H}_2$. Suppose $K\in\M_{\ell,m}$. Then there exists $C>0$,
  $\epsilon_0>0$, such that
  \[
  \norm{K(\epsilon,h)}\leq C\epsilon^\ell h^m, \qquad \forall
  \epsilon\in(0,\epsilon_0], \quad \forall
  h^{-1}\in\frac{\abs{\nu}}{2}+\pscal{\nu}{\NM^n} .
  \]
\end{lemm}

To prove this lemma we shall use a specific version of the functional
calculus for the harmonic oscillator, which can be obtained as a small
improvement of the semiclassical functional calculus
of~\cite{dimassi-sjostrand}. The only difference is that we allow any
value of $h_0$, whereas in~\cite{dimassi-sjostrand} $h_0$ has to be
chosen ``small enough''. This modification is essential for our
purposes.
\begin{lemm}
  \label{lemm:fonctionnel}
  Let $\hat{H}_2$ be a harmonic oscillator as
  in~\eqref{equ:harmonic_osc}. Then for any $h_0>0$ and for any
  $f\in\Cinf_0(\RM)$, the family of operators $(f(\hat{H}_2))_{h\leq
    h_0}$ is a semiclassical pseudo-differential operator in $\Psi(0)$:
  there exists a bounded family $(a(\cdot,h))_{h\leq h_0}$ in $S(0)$
  such that $f(\hat{H}_2)=Op_h^w(a(\cdot,h))$. Furthermore, for any
  integer $m$,
  \[
  a(\cdot,h)=a_0+h a_1+h^2a_2+\cdots+h^mr_m(\cdot,h)
  \]
  where all $a_i$ are $\Cinf_0$ with support in the support of
  $f(H_2)$ (and $a_0=f(H)$, $a_1=0$), and $(r_k(\cdot,h))_{h\leq h_0}$
  is a bounded family in $S(-N)$, for all $N\geq 0$.
\end{lemm}

\begin{demo}
  Assume the result is true for $h_0'$ and let us prove it for $h_0 =
  \delta^2 h'_0$. Let $U_\delta$ be the unitary map $(U_\delta \Psi )
  (x) = \delta^{\frac{n}{2}} \Psi (\delta x)$. Since
$$ U_\delta \hat{H}_2( h ) U_\delta^* = \delta^2 \hat{H}_2( h/ \delta^2 ) $$
and $U_\delta$ is unitary, we have
\[
U_\delta f( \hat{H}_2 ( h )) U_\delta^* = f( \delta^2 \hat{H}_2(
h / \delta^2 ))
\]
By the known result, $f( \delta^2 \hat{H}_2( h ')) = Op^w_{h'} (
b(.,h'))$ for $h' \leqslant h_0'$. So
$$ f( \hat{H}_2 (h )) = Op^w_h (a(. ,h) ), \quad \text{ for } h
\leqslant h_0 $$ with $a(x,\xi, h) = b( \delta^{-1} x,\delta^{-1} \xi,
\delta^{-2} h)$. And we deduce the estimates of $a$ from those of $b$.
\end{demo}

\begin{demo}[of lemma~\ref{lemm:normes}]
  Let $k(x,\xi;\h)$ be the Weyl symbol of $K$. First (as usual) we can
  assume that $k$ has compact support. Indeed we split $k=k^b+k^\infty$,
  where $k^b$ has compact support, $k^\infty$ vanishes on
  \[
  B_{E_0+\delta}:=\{z; H_2(z)\leq E_0+\delta\}
  \]
  for some positive $\delta$, and both $k^b$ and $k^\infty$ commute
  with $H_2$ (using a $\T^k$-averaging as in the proof of
  proposition~\ref{prop:weyl}). Then we choose $f\in\Cinf_0(\RM)$ with
  value $1$ on $[0,E_0]$ and support inside $(-\infty,E_0+\delta)$, so
  that $f(\hat{H}_2)$ and $K^\infty:=Op_\h^w(k^\infty)$ have disjoint
  microsupport. Then, by lemma~\ref{lemm:fonctionnel} and symbolic
  calculus, $f(\hat{H}_2)K^\infty=\ohb$ uniformly for $\h\leq 1$.
  Since $\h=h\epsilon^2$, $K^\infty(\epsilon,h)$ is of order
  $\O(h^N\epsilon^N)$ for all $N$, provided $h\epsilon^2\leq
  2E_0/\abs{\nu}$.
  Now consider $q(z;\epsilon,h):=k^b(\epsilon z;\epsilon^2 h)$. Since
  $k^b$ has compact support and its Taylor expansion is in $\O_\ell$,
  we have
  \[
  q = \sum_{2\alpha+\beta=\ell}\O((\epsilon^2 h)^\alpha\abs{\epsilon
    z}^\beta) = \epsilon^\ell
  \sum_{2\alpha+\beta=\ell}\O(h^\alpha\abs{z}^\beta).
  \]
  Since $k^b=\O(\h^m)$ we know that we can actually write
  \[
  q = \epsilon^\ell h^m
  \sum_{2\alpha'+\beta=\ell-2m}\O(h^{\alpha'}\abs{z}^\beta).
  \]

  We could conclude here by restricting to a compact subset in the new
  variable $z$ and using standard semiclassical arguments similar to
  what we already used at several places. However this would require
  to restrict the validity of the expected estimate to sufficiently
  small values of $h$. So we need to refine the analysis.

  By hypothesis $h$ is bounded (by $2/\abs{\nu}$); so there is a
  constant $C>0$ (depending of course on $\ell,m,\nu$) such that
  \[
  \abs{q(z;\epsilon,h)}\leq C\epsilon^\ell h^m\bracket{z}^{\ell-2m} ,
  \]

  For all multi-indices $\alpha\in\NM^{2n}$, $\partial^\alpha_z k^b$
  still has compact support and its Taylor expansion is in
  $\O_{\ell-\abs{\alpha}}$. So we see similarly that for any $\alpha$,
  \[
  \abs{\partial^\alpha_z q(z;\epsilon,h)} \leq C_\alpha\epsilon^\ell
  h^m\bracket{z}^{\ell-2m-\abs{\alpha}} .
  \]
  In particular $(h^{-m}\epsilon^{-\ell}Op_h^w(q))\in
  S(\ell-2m)$, uniformly for all $\epsilon\leq
  \epsilon_0$, $h\leq 2/\abs{\nu}$.  Using
  lemma~\ref{lemm:fonctionnel} with $h_0=2/\abs{\nu}$ and a function
  $f\in\Cinf_0(\RM)$, we consider the operator
  \begin{equation}
    \tilde{Q}:=h^{-m}\epsilon^{-\ell} Op^w_h(q) f(Op^w_h(H_2)).
    \label{equ:Q1}
  \end{equation}
  
  By the composition theorem for pseudo-differential operators
  (\cite[proposition 7.7]{dimassi-sjostrand}), $\tilde{Q}$ is
  uniformly in $S(0)$, and hence bounded.

  It remains to recall that $K(\epsilon,h)=Op^w_h(q)_{\restr \H_h}$
  which says that
  \begin{equation}
    K(\epsilon,h)=h^m\epsilon^\ell f(1)\tilde{Q},
    \label{equ:Q2}
  \end{equation}
  and gives the result.
\end{demo}

Let us now turn to some applications of our results.
\subsection{Polynomial differential operators}
If $k\in\D_i$ then $Op^w_\h(k)$ is a $\h$-differential operator with
polynomial coefficients. Using~\eqref{equ:scaling-weyl} we see that
\begin{equation}
  U_\epsilon Op^w_\h(k)U_\epsilon^* = \epsilon^i Q(h),
  \label{equ:scaling-epsilon}
\end{equation}
where $Q(h)$ is a $h$-differential operator with polynomial
coefficients, independent of $\epsilon$, of order at most $i$ in
$(x,\partial_x)$.

Now suppose $K$ is a $\h$-pseudo-differential operator commuting with
$\hat{H}_2$, with formal Weyl symbol $[K]$. Fix $\ell\geq 0$ and write
\[ [K] = k^{(0)}+k^{(1)}+\cdots+k^{(\ell)} + r_{\ell+1} ,
\]
with $k^{(i)}\in\D_i$ and $r_{\ell+1}\in\O(\ell+1)$. Applying
lemma~\ref{lemm:normes}, we get
\begin{equation}
  K(\epsilon,h) = K_0(h)+\epsilon K_1(h) + \cdots + \epsilon^\ell
  K_\ell(h) + \epsilon^{\ell+1}R_{\ell+1}(h),
  \label{equ:K}
\end{equation}
where $R_{\ell+1}(h)$ is uniformly bounded with respect to $h$ and
$K_i(h)=Op^w_h(k^{(i)})$ (actually, with some abuse of notation, we
have written $K_i(h)$ for their restrictions to $\H_h$).

In other words, in order to study the spectrum of $K(\epsilon,h)$, we
can restrict to the study of polynomial differential operators
commuting with $H_2$, up to an error of order $\O(\epsilon^{\ell+1})$,
uniformly in $h$. This is very important in practice for numerical
calculations, since the matrix elements of such $K_i$'s on a basis of
eigenvectors of the harmonic oscillator are very easy to compute
(especially when the $K_i$'s are written in terms of creation an
annihilation operators or --- this amounts to the same --- when
studied in the Bargmann-Fock representation). See \cite{san-colin} for
the case of the $1:2$ resonance.

We state the result with the original parameters $(E,\h)$ for further
reference.
\begin{theo}
  \label{theo:FNB-poly}
  Let $P$ be a $\h$-pseudo-differential operator satisfying the
  hypothesis of theorem~\ref{theo:FNB}. Then there exists a harmonic
  oscillator $H_2$ on $\RM^{2n}$ and a formal series
  \[ [K]=k^{(3)}+k^{(4)}+\cdots \in\O_3
  \]
  commuting with $H_2$ such that, for any $\ell\in\NM$ and any
  $\eta>0$, there exists $E_0>0$, $\h_0>0$, $C>0$ such that for all
  $E\in(0,E_0]$, $\h\in(0,\h_0]$ and $j\in\NM$,
  \[
  \left(\lambda^P_j\leq E \textrm{ or } \lambda^{Q_{\ell}}_j\leq E
  \right) \impliq \abs{\lambda^P_j-\lambda^{Q_\ell}_j - \h H_0}\leq
  CE^{\frac{\ell+1}{2}} + \O(\h^\infty),
  \]
  where
  \[
  Q_{\ell} :=
  (\hat{H}_2+Op^w_\h(k^{(3)}+k^{(4)}+\cdots+k^{(\ell)}))_{\restr
    \Pi^{\hat{H}_2}_{(-\infty, (1+\eta)E]}(L^2(\RM^n))} .
  \]
\end{theo}
\begin{demo}
  Naturally, we first apply theorem~\ref{theo:FNB}.  Therefore the
  question is reduced to the determination of the spectrum of the
  operator $Q((1+\eta)E)$ defined in that theorem, up to an error of
  order $\O(E^\infty)+\O(\h^\infty)$. We need to compare
  $Q((1+\eta)E)$ to $Q^{\ell}$. But equation~(\ref{equ:K}) yields
  \[
  Q((1+\eta)E)_{\restr \H_{E,\h}}-(Q_{\ell})_{\restr \H_{E,\h}} =
  \epsilon^{\ell+1}U_\epsilon^*R_{\ell+1}(h)U_\epsilon =
  \O(\epsilon^{\ell+1})
  \]
  with $\epsilon=\sqrt{E}$. Hence $Q((1+\eta)E)-Q_{\ell} =
  \O((1+\eta)E)^{\frac{\ell+1}2}=\O(E^{\frac{\ell+1}2})$. The minimax
  (lemma~\ref{lemm:comparaison}) gives the result.
\end{demo}

\subsection{The semi-excited Weyl law}
Let $P$ be a $\h$-pseudo-differential operator satisfying the
hypothesis of theorem~\ref{theo:FNB}, and let $H_2$ be the
corresponding harmonic oscillator, as in~\eqref{equ:harmonic_osc}. We
define the \emph{resonance order} of $H_2$ to be
\[
r_\nu:=\max(3,\inf\{\abs{\alpha}; \quad \alpha\in\ZM^n, \alpha\neq 0,
\pscal{\alpha}{\nu}=0\})\in\NM\cup\{\infty\}.
\]
By convention $r_\nu=\infty$ is $\nu$ is non resonant. Notice that, in
contrast to corollary~\ref{coro:resonances}, we exclude here the
value $r_\nu=2$ since we shall always deal with perturbations terms
that are formally in $\O_3$. This remark pertains through the
remaining of the article.

We are interested here in the counting function
\[
\N^P(E,\h):=\#\{j; \quad \lambda^P_j(\h)\leq E\}.
\]
\begin{theo}
  \label{theo:weyl}
  Let $p$ be the principal symbol of $P$. For any $\ell<r_\nu$ (and
  $\ell\geq 2$)
  \[
  \N^P(E,\h)=\frac{1}{(2\pi\h)^n}\int_{p\leq E}\abs{dx d\xi} +
  \O(\h^{-n}E^{n-1}(\h+E^{\frac{\ell+1}{2}}))
  \]
  uniformly for $E$ and $\h$ small enough.
\end{theo}
\begin{coro}
  \label{coro:weyl}
  For any $C>0$,$\gamma>0$, $\ell\in[2,r_\nu)$ the Weyl counting
  function for ``semi-excited states'' has the following asymptotics,
  as $\h\fleche 0$:
  \[
  \N^P(C\h^\gamma;\h) = \frac{1}{(2\pi\h)^n}\int_{p\leq E}\abs{dx
    d\xi} + \O(\h^{(1-n)(1-\gamma)}(1+\h^{\gamma\frac{\ell+1}{2}-1}))
  .
  \]
\end{coro}
In other words if $\gamma\leq \frac{2}{\ell+1}$ then the remainder is
$\O(\h^{-n+\gamma(n+\frac{\ell-1}{2})})$ while if $\gamma\geq
\frac{2}{\ell+1}$ the remainder is $\O(\h^{(1-n)(1-\gamma)})$.  We see
that if $\gamma<\frac{2}{r_\nu+1}$ then one always gets the best
remainder $\O(\h^{(1-n)(1-\gamma)})$.

\begin{demo}[of theorem~\ref{theo:weyl}]
  Apply theorem~\ref{theo:FNB-poly} with the same $\ell$. Then, using
  the notation of that theorem,
  \begin{equation}
    \N^P(E,\h) = \N^{Q_\ell}(E,\h) + r(E,\h),
    \label{equ:weyl-reste}
  \end{equation}

  where $r(E,\h)\leq\#\{j;\quad
  \lambda^{Q_\ell}_j\in[E-CE^{\frac{\ell+1}{2}},
  E+CE^{\frac{\ell+1}{2}}]\}$.  Since $\ell<r_\nu$ we can see using
  corollary~\ref{coro:resonances} that
  \[
  Q_\ell =
  f\left(\frac{x_1^2+\hat{\xi}_1^2}{2},\dots,\frac{x_n^2+\hat{\xi}_n^2}{2};\h\right)
  ,
  \]
  where $f(u;\h)$ is a polynomial in $(u,\h)$ of degree at most
  $[\ell/2]$, with linear part equal to $\pscal{u}{\nu}$. The
  eigenvalues of $Q_\ell$ are therefore of the form
  \[
  \lambda_j^{Q_\ell} =
  f\left(\h(\sfrac12+\alpha_1),\dots,\h(\sfrac12+\alpha_n)\right)
  \]
  for integers $(\alpha_1,\dots,\alpha_n)\in\NM^n$, and
  \begin{eqnarray*}
    &\N^{Q_\ell}(E_1,E_2;\h)  := \#\{j;\quad
    \lambda^{Q_\ell}_j\in[E_1,E_2]\}&\\
    &=  \#\{\alpha\in\NM^n;\quad
    f\left(\h(\sfrac12+\alpha_1),\dots,\h(\sfrac12+\alpha_n)\right) \in
    [E_1,E_2]\}\nonumber .&
    \label{equ:weyl}
  \end{eqnarray*}
  Using a covering of $\RM^n$ by small cubes with sides of length
  $2\h$, centred at
  $\left(\h(\sfrac12+\alpha_1),\dots,\h(\sfrac12+\alpha_n)\right)$,
  $\alpha\in\ZM^n$, it is then easy to see that
  \begin{eqnarray}
    & \h^n\N^{Q_\ell}(E_1,E_2;\h) =
    \textup{Vol}(f^{-1}([E_1,E_2]))&\nonumber \\
    \label{equ:weyl-vol}
    & + \O(\textup{Vol}(f^{-1}([E_1-\h,E_1+\h]))) +
    \O(\textup{Vol}(f^{-1}([E_2-\h,E_2+\h])))  &
  \end{eqnarray}
  For $E$ small enough, using that $f(u;\h)\sim \pscal{u}{\nu}$ with
  $\nu_i>0$, one has the bound
  \[
  \textup{Vol}(f^{-1}([E-\delta,E+\delta]))=\O(\delta E^{n-1}) .
  \]
  Applying this with $\delta=\h$ and $\delta=CE^{\frac{\ell+1}{2}}$ we
  get from~(\ref{equ:weyl-vol})
  \[
  \N^{Q_\ell}(E-CE^{\frac{\ell+1}{2}}, E+CE^{\frac{\ell+1}{2}}) =
  \O(\h^{-n}E^{n-1}(\h+E^{\frac{\ell+1}{2}})).
  \]
  Another application of~(\ref{equ:weyl-vol}), with $E_1=0$, $E_2=E$,
  combined with~\eqref{equ:weyl-reste}, yields:
  \[
  \N^P(E;\h)=\h^{-n}\textup{Vol}\{f\leq E\} +
  \O(\h^{-n}E^{n-1}(\h+E^{\frac{\ell+1}{2}})).
  \]
  Notice now that, by a simple Fubini argument,
  \[
  \textup{Vol}\{f\leq E\}=(2\pi)^{-n}\int_{q^{\ell}\leq E}\abs{dx
    d\xi},
  \]
  where $q^{\ell}$ is the full symbol of $Q_\ell$. We know from
  lemma~\ref{lemm:Rn} and proposition~\ref{prop:weyl} that
  $p=q^{\ell}+\O_{\ell+1}+\O(\h)$.  Therefore
  \[
  \{q^{\ell}\leq E - C'(\h+E^{\frac{\ell+1}{2}})\}\subset \{p\leq
  E\}\subset \{q^{\ell}\leq E + C'(\h+E^{\frac{\ell+1}{2}})\},
  \]
  so $\int_{q^{\ell}\leq E}\abs{dx d\xi} = \int_{p\leq E}\abs{dx d\xi}
  + \O(E^{n-1}(\h+E^{\frac{\ell-1}{2}}))$.
\end{demo}

\subsection{The low-lying eigenvalues}
Let $P$ be a $\h$-pseudo-differential operator satisfying the
hypothesis of theorem~\ref{theo:FNB}, and let $H_2$ be the
corresponding harmonic oscillator, as in~\eqref{equ:harmonic_osc}.
For an energy $E$ of order $\h$, the Weyl formula
(corollary~\ref{coro:weyl}) says that the number of eigenvalues of $P$
below $E$ is bounded, independently of $\h$. Of course, this can be
obtained directly as a consequence of theorem~\ref{theo:FNB-poly}:
because for $E=C\h$, the dimension of $\Pi^{\hat{H}_2}_{(-\infty,
  (1+\eta)E]}(L^2(\RM^n))$ is independent of $\h$.  We use here the
Birkhoff normal form to recover a result of
Helffer-Sj{\"o}strand~\cite[theorem 3.6]{helffer-sjostrand} concerning
the asymptotics of the low-lying eigenvalues of Schr{\"o}dinger operators.
\begin{theo}
  For any $C>0$ and $\h$ small enough, the spectrum of $P$ in
  $(-\infty,C\h]$ consists of a finite number, independent of $\h$, of
  eigenvalues. These eigenvalues admit an asymptotic expansion of the
  form

  \begin{equation}
    \lambda_j^P(\h)\sim \h H_0 + \h\mu_0 + \h^{\frac{3}{2}}\sum_{m=0}^\infty
    \h^{\frac{m}{2}}\mu_m,
    \label{equ:fond-de-puits}
  \end{equation}
  where $\h\mu_0$ is an eigenvalue of $\hat{H}_2$, and, as usual,
  $H_0$ is the value of the sub-principal symbol of $P$ at the minimum
  of its principal symbol.  The number of eigenvalues with a given
  $\mu_0$ is equal for $\h$ small enough to the multiplicity of
  $\h\mu_0$ for $\hat{H}_2$. In particular the smallest eigenvalue of
  $P$ is simple for $\h$ small enough.
\end{theo}
\begin{demo}
  We apply theorem~\ref{theo:FNB-poly} with $E=C\h$ and $\ell$ large.
  Thus, modulo an error of size $\O(\h^{(\ell+1)/2})$, the eigenvalues
  of $P$ less than $C\h$ are equal to the eigenvalues of $Q_\ell+\h
  H_0$.  Decomposing the Hilbert space on which $Q_\ell$ acts
  according to the eigenspaces of the harmonic oscillator $\hat{H}_2$,
  and using the unitary equivalence as in~\eqref{equ:scaling-epsilon},
  we are reduced to the study of eigenvalues of matrices of the form
  \[
  Q^\ell(\epsilon,h) = \epsilon^2 K_2(h)+\epsilon^3 K_3(h) + \cdots +
  \epsilon^\ell K_\ell(h) ,
  \]
  acting on $\H_h$, where
  $h^{-1}=\frac{\abs{\nu}}{2}+\pscal{\alpha}{\nu}$ for some fixed
  $\alpha\in\NM^n$ and $\epsilon=\sqrt{\h h^{-1}}$. By standard
  perturbation theory for matrices, the spectrum of
  $\epsilon^{-2}Q^{\ell}(\epsilon,h)$ is analytic in $\epsilon$ for
  small $\epsilon$.  Since $\epsilon^2 K_2(h)$ is unitarily equivalent
  to $\hat{H}_2$ we obtain the expansion~\eqref{equ:fond-de-puits},
  along with the statement concerning the multiplicities.
\end{demo}

\begin{rema}
  In case $P$ is a semiclassical Schr{\"o}dinger operator
  $-\frac{\h^2}{2}\Delta+V(x)$, This result appeared almost
  simultaneously in~\cite[theorem 5.1]{simon1} and ~\cite[theorem
  3.6]{helffer-sjostrand}.  The techniques of~\cite{helffer-sjostrand}
  could in principle be easily generalised to treat, like we do here,
  general pseudo-differential operators. While these results are now
  well known, is it interesting to remark that the appearance of
  \emph{half-integer} powers of $\h$ in the asymptotic expansion was
  not so obvious at that time. Actually in~\cite{simon1} these
  half-integer powers had been forgotten.

  The Birkhoff normal form we used here makes it very clear as to why
  and when such powers may appear in the asymptotic expansion of the
  eigenvalues.  In particular the smallest exponent from which such
  half-integer powers can appear is half the resonance order $r_\nu$.
  Indeed, below this order, corollary~\ref{coro:resonances} shows that
  the operators $K_j(h)$ must be of the form
  $f_j(\frac{x_1^2+\hat{\xi}_1^2}{2},\cdots,\frac{x_n^2+\hat{\xi}_n^2}{2})$
  for some polynomial $f_j$ and hence have even order in $(x,\xi)$.
  Notice that in dimension 1, $r_\nu=\infty$; hence only integer
  powers of $\h$ may show up in that case. As remarked
  in~\cite{helffer-sjostrand}, the simplest case where half-integer
  powers of $\h$ are present is the so-called $1:2$ resonance: $n=2$
  and $\nu=(1,2)$. The coefficient of $\h^{3/2}$ is then the average
  along the flow of $H_2$ of the term of order 3 in the Taylor
  expansion of the symbol. A more general statement is given in
  section~\ref{sec:asymptotics}.
\end{rema}

\section{Toeplitz operators}
\label{sec:toeplitz}

As we saw in equation~\eqref{equ:torus}, our initial spectral problem
is reduced to the spectral analysis of a pseudo-differential operator
invariant under a $\T^k$ action, for some $k=1,\dots,n$. With this
respect, we have two qualitatively extreme situations: $k=n$ and
$k=1$. In the first case $k=n$ the harmonic oscillator $H_2$ has no
resonance relation, and the situation is essentially completely
integrable, at least in the semi-excited regime. This can be seen
explicitly with theorem~\ref{theo:FNB-poly}, in view of
corollary~\ref{coro:resonances}.

On the contrary, in the case $k=1$, the harmonic oscillator $H_2$ is
completely resonant: up to a common multiple, all frequencies $\nu_i$
are integers. The $\T^1$ action is precisely the flow of $H_2$. From
the perspective of integrability, the spectral analysis looks more
involved. However the fact that $H_2$ has a periodic flow forces the
spectrum to split into regularly spaced clusters. The spectral
analysis becomes in some sense simpler, since it comes down to the
study of each individual cluster. From a geometric point of view each
of these clusters correspond to the spectrum of an operator acting on
a reduced space. The goal of this section is to develop this idea.

In principle, the mixed case $1<k<n$ could be treated by a combination
of both techniques. But this still has to be investigated further.

\ouf

We use the notation introduced in the beginning of
section~\ref{sec:joint}. In particular we deal with rescaled
parameters $\epsilon,h$ and, using the scaling operator $U_\epsilon$
as in~\eqref{equ:scaling}, we are able to reduce the analysis to the
case $E=1$. The new tool we introduce here is to work in the Bargmann
space $\B_h$~\cite{bargmann}. Recall that $\B_h$ is the space of
entire holomorphic functions on $\CM^n$ with finite norm, where the
norm comes from the scalar product
\[
\pscal{\psi}{\psi'}_{\CM^n} = \int_{\CM^n}(\psi,\psi')(z)\mu(z),
\qquad \textrm{ with }
(\psi,\psi')(z)=\psi(z)\overline{\psi'(z)}e^{-\abs{z}^2/h}
\]
where $\abs{z}^2=\sum z_i\bar{z_i}$ and $\mu$ is the Lebesgue measure
on $\CM^n=\RM^{2n}$. Operators on $L^2(\RM ^n)$ can be transported on
$\B_h$ via the Bargmann transform which is the unitary map
$U_{\B}:L^2(\RM^n)\fleche\B_h$ given by
\begin{gather}
  \label{equ:transformee-Bargmann}
  U_{\B}(\phy)(z) = \frac{2^{n/4}}{(2\pi
    h)^{3n/4}}\int_{\RM^n}e^{h^{-1}\left( z \cdot x \sqrt{2} -
      (z^2+x^2)/2 \right)}\phy(x)dx ,
\end{gather}
where $z\cdot x=\sum z_ix_i$, $z^2=z\cdot z$, $x^2=x\cdot x$.

The harmonic oscillator $\hat{H}_2 (h) :=Op^w_h(H_2)$ with $H_2$ as
in~\eqref{equ:harmonic_osc} becomes
\[
\hat{H}_2^{\B} (h) :=U_{\B}\tilde{H}_2 (h) U_{\B}^*=h\sum_{j=1}^n
\nu_j\left(z_j\deriv{}{z_j}+\sfrac{1}{2}\right) .
\]
In order to deal with symbols of operators in the Bargmann side we
simply identify the real phase space $\RM^{2n}$ with $\CM^n$ using
$z_j=(x_j-i\xi_j)/\sqrt{2}$.
Hence $H_2=\sum_j \nu_j\abs{z_j}^2$.

\subsection{The reduction setting}

In all this section~\ref{sec:toeplitz} the main assumption is that for
all $i,j$, $\nu_i/\nu_j\in\QM$.  Then there exists positive coprime
integers $\p_1,\dots,\p_n$, and a constant $\nu_c>0$ such that
$\nu_j=\nu_c\p_j$. Hence, according to~\eqref{equ:rescaling}, our
rescaled semiclassical parameter $h$ is of the form
\[
h^{-1} = \frac{\abs{\nu}}{2} + \nu_cN, \quad N\in\NM .
\]

Let $Y=\{y\in\CM^n; \quad H_2(y)=1\}$; it is a smooth, compact
submanifold of $\CM^n$.  The Hamiltonian flow of $H_2$ defines a
locally free action of $S^1$ on $Y$
\begin{gather} \label{eq:action} S^1\times Y \fleche Y, \qquad
  (u,y)\mapsto u.y:=(u^{\p_1}y_1,\dots,u^{\p_n}y_n),
\end{gather}
where we identify $S^1$ with the complex numbers of modulus $1$.

The quotient $Y/S^1$ is an orbifold $M$. It is endowed with the
Marsden-Weinstein symplectic form $\omega$, naturally defined by
$\pi^*\omega = \omega_Y$ where $\pi$ is the projection $Y\fleche M$
and $\omega_Y$ is the restriction to $Y$ of the canonical symplectic
form of $\CM^n$. Smooth functions on $M$ are by definition
$S^1$-invariant functions on $Y$.

Let $\Cinf_N(Y)\subset \Cinf(Y)$ be the space of equivariant functions
in the following sense:
\[
\psi(u.y) = u^N\psi(y), \quad \forall u\in S^1, y\in Y.
\]

$M$ is naturally endowed with a complex line bundle $L\fleche Y$ whose
sections are identified with functions of $\Cinf_1(Y)$ using the
pull-back $\pi^*$. More generally $\Cinf(M,L^N)\simeq \Cinf_N(Y)$.
Let us endow $L$ with the hermitian structure such that, if $\psi_r$
and $\psi_r'$ are sections of $L^N$ and $\psi=\pi^*\psi_r$,
$\psi'=\pi^*\psi_r'$, then
\begin{gather} \label{eq:produit_scalaire} (\psi_r,\psi_r')(\pi(y)) =
  \psi(y)\overline{\psi'(y)}e^{-\nu_c N\abs{y}^2}.
\end{gather}
This defines the scalar product $\pscal{\psi_r}{\psi'_r}_M:=\int_M
(\psi_r,\psi_r') \mu_M$, where $\mu_M$ is the Liouville (or
symplectic) measure of $M$.

On the other hand let $\H_h^{\B} = \ker (\hat{H}^\B_2 (h) -1) = U_\B
\H_h$ where $\H_h$ is, as in section~\ref{sec:joint}, the eigenspace
of $\hat{H}_2 (h)$ corresponding to the eigenvalue $E=1$. It is well
known~\cite{bargmann} that the monomials
$z^\alpha=z_1^{\alpha_1}\cdots z_n^{\alpha_n}$ such that
\begin{equation}
  \pscal{\p}{\alpha} = N
  \label{equ:monomes}  
\end{equation}
form a basis of $\H_h^{\B}$. But equation~\eqref{equ:monomes} also
shows that the restrictions of the $z^\alpha$'s to $Y$ belong to
$\Cinf_N(Y)$. Projecting onto $M$ we thus see that $\H_h^\B$ may be
identified with a subspace $\BH_N\subset \Cinf(M,L^N)$. Since in
general this bijection
\[
V_N : \H^\B_h \fleche \BH_N
\]
is not unitary~\cite{charles-reduc}
we introduce $U_N=(V_N V_N^*)^{-\frac12}V_N : \H_h^\B\fleche\BH_N$.
The sequence of spaces $(\BH_N)_{N\in\NM}$ must be viewed as the
quantising Hilbert space for the reduced phase space $M$, and $U_N$ as
a kind of Fourier integral operator allowing to transport equivariant
wave functions of the original Bargmann space onto the reduced Hilbert
space $\BH_N$.

\subsection{Reduction of the Birkhoff normal form}

Let $K=K(\h)$ be an $\h$-pseudo-differential operator on $\RM^n$
commuting with $\hat{H}_2$. Assume $K$ belongs to the class
$\M_{\ell,m}$, as in definition~\ref{defi:Mlm}, and consider again the
family of rescaled operators $K(\epsilon,h) : \H_h\fleche \H_h$. Since
the Taylor expansion of the Weyl symbol $[K]$ is $\O_\ell$ we see as
in~\eqref{equ:K} that the family $\epsilon^{-\ell}K(\epsilon,h)$ is
continuous at $\epsilon=0$, if we set
\[
\epsilon^{-\ell} K(\epsilon,h)_{\restr \epsilon=0} = K_\ell(h),
\]
where $K_\ell(h)$ is the differential operator with polynomial
coefficients obtained by the leading term of $[K]$, as precisely
defined in~\eqref{equ:K}.

We wish to consider $K(\epsilon,h)$ as a reduced operator on $M$. For
this purpose, let $K^\B(\epsilon,h)=U_\B K(\epsilon,h) U_\B^*$ and
introduce
\[
T(\epsilon,h) := \epsilon^{-\ell} U_N K^\B(\epsilon,h) U_N^* :
\BH_N\fleche \BH_N.
\]

Let $\BPi_N$ denote the orthogonal projector on $\BH_N$.

The main result of this section is that $T$ is a semiclassical
Toeplitz operator on $M$, as stated in the following theorem.
\begin{theo} \label{theo:toeplitz} If $K\in\M_{\ell,m}$, then there
  exists a sequence of functions $(f(\cdot;N))_{N\in\NM}$ in
  $\Cinf(M\times[0,\epsilon_0])$ admitting an asymptotic expansion
  \[
  f(x,\epsilon;N) = f_0(x,\epsilon) + N^{-1}f_1(x,\epsilon) +
  N^{-2}f_2(x,\epsilon) + \cdots
  \]
  for the $\Cinf$ topology, such that
  \[
  T(\epsilon,h) = N^{-m}\BPi_N f(\cdot,\epsilon;N) + \O(N^{-\infty}),
  \]
  uniformly in $\epsilon$. Moreover, if $k_m$ is the $\h$-principal
  symbol of $K$ then
  \[
  \epsilon^\ell f_0(\pi(z),\epsilon) = k_m(\epsilon z), \qquad \forall
  z\in Y.
  \]
\end{theo}

\begin{demo}
  The proof proceeds by first representing $K$ as a Toeplitz operator
  on $\CM^n$ using the Bargmann transform and the corresponding
  formula for the Toeplitz symbol. Then one has to show that this
  symbol can be restricted to the energy hypersurface $H_2=1$ while
  retaining the asymptotics in $h$ (or $N^{-1}$) and the $S^1$
  invariance.

  However, some technical preliminaries are required before this.  As
  in the proof of lemma~\ref{lemm:normes} we can assume that the
  symbol $k=k(x,\xi;\h)$ is compactly supported. Again we introduce
  the rescaled symbol
$$q(x,\xi;\epsilon,h)=k(\epsilon
x,\epsilon\xi;\epsilon^2 h).$$ Refining the proof of
lemma~\ref{lemm:normes} we can check that $q(\cdot,\epsilon,h)$ can be
assumed to have a support $D$ that does not depend on $\epsilon$ and
$h$ and with $0\not\in D$. Indeed, let $\Phi\in\Cinf_0(\RM^{2n})$ with
support not containing the origin and identically equal to 1 on a a neighbourhood of $\{H_2=1\}$.
By $h$-symbolic calculus
  $$(1-Op^w_h(\Phi))\phy(Op^w_h(H_2))= \O(h^\infty),$$ 
  provided the support of $\phy$ is suitably restricted around $1$.
  Hence if we replace $q$ by the Weyl symbol of
  $Op^w_h(q)Op^w_h(\Phi)$, we deduce from equations~\eqref{equ:Q1}
  and~\eqref{equ:Q2} that $K(\epsilon,h)$ is modified by a term of
  order $\O(\epsilon^\ell h^\infty)$, which proves our claim. Moreover
  one has
  \[
  q(x,\xi;\epsilon,h) = h^m\epsilon^\ell \tilde{q}(x,\xi;\epsilon,h)
  \]
  where $\tilde{q}(\cdot,h)$ is a family of
  $\Cinf_0(\RM^{2n}\times[0,\epsilon_0])$ admitting an asymptotic
  expansion in powers of $h$ for the $\Cinf$ topology. Hence, for the
  proof of the theorem, we can safely assume $m=\ell=0$.

  We can now turn to the Bargmann side. Let
  \[
  Q^\B = U_\B Op^w_h (q(\cdot;\epsilon,h)) U_\B^*.
  \]
  We know that $Q^\B$ can be represented as a Toeplitz operator.
  Precisely, let $\BPi^\B$ be the orthogonal projector of
  $L^2(\CM^n,e^{-\abs{z}^2/h})$ onto $\B_h$. For any bounded function
  $g$ on $\CM^n$, the Toeplitz operator with contravariant symbol $g$
  is by definition the operator
  \[
  T_g : \B_h \fleche \B_h, \qquad \psi \mapsto \BPi^\B (g\psi).
  \]
  Then we have the following result
  \begin{theo} \label{sec:Weyl-contravariant} If $g$ is in the symbol
    class $S(0)$ on $\CM^n\simeq \RM^{2n}$ then the operator $U^*_\B
    T_g U_\B : L^2(\RM^n)\fleche L^2(\RM^n)$ is a pseudo-differential
    operator whose Weyl symbol is
    \[
    I(g)(\zeta) = \frac{1}{(\pi h)^n}\int_{\CM^n}
    e^{-2h^{-1}\abs{z}^2} g(\zeta+z)\abs{dz d\bar{z}}.
    \]
    The map $I$ is continuous from $S(0)$ to $S(0)$. Moreover for any
    $g\in S(0)$ and all $k\geq 1$,
    \begin{equation}
      I(g) = \sum_{j=0}^{k-1} (\sfrac{h}{2})^j \Delta^jg/j! + h^kR_k(g)
      \label{equ:toeplitz-weyl}
    \end{equation}
    where $R_k$ is a continuous map from $S(0)$ to $S(0)$.
  \end{theo}
  \begin{demo}
    Assume first that the symbol $g$ is in the Schwartz class. Since
    $U^*_\B U_\B = \Id $ and $U_\B U^*_\B = \BPi^\B$, one has $U^*_\B
    T_g U_\B = U^*_\B g U_\B $ and its kernel is
$$ K (x,y) = \int_{\CM^n} U^*_\B (x,z)  g (z)  U_\B (z,y) \mu (z) $$
The kernel $U_\B (z,y)$ is given in (\ref{equ:transformee-Bargmann}),
$U^*_\B (x,z) = e^{-h^{-1} \abs{z}^2} \overline{U_\B (z,x)}$ so
$$ K (x,y) = \frac{2^{n/2}}{(2 \pi h)^{3n/2}} \int_{\RM^{2n}} e^{i h^{-1}
  \left(p \cdot (x-y) +\frac{i}{4} (x-y)^2 +i (q -\frac{1}{2}(x+y))^2
  \right) }g (q,p) |dq dp| $$ with $z = (q - ip ) /\sqrt{2}$. We
recover the Weyl symbol $I(g)$ with the well-known formula
\begin{xalignat*}{2}
  I(g)(x, \xi) =& \int_{\RM^n} e^{-i h^{-1} v.\xi} K \Bigl( x+
  \frac{v}{2}
  , x - \frac{v}{2} \Bigr) |dv| \\
  = & \frac{2^{n/2}}{(2\pi h)^{3n/2}} \int_{\RM^{3n}} e^{i h^{-1}
    \left( v.(p - \xi) + i ( q- x ) ^2 + \frac{i}{4} v^2 \right)} g
  (q,p) |dq
  dpdv| \\
  = & \frac{2^{n/2}}{(2 \pi h)^{3n/2}} \int_{ \RM^{3n} } e^{i h^{-1}
    \left( v.p + \frac{i}{4} v^2 + i q^2 \right) } g ( x+ p , \xi + q)
  |dq
  dp dv|  \\
  = & \frac{1}{(\pi h)^n}\int_{\RM^{2n}} e^{-h^{-1} \left( p^2 + q^2
    \right)} g ( x+ p , \xi + q) |dq dp|
\end{xalignat*}
Going back to the coordinates $z$ and $\zeta = (x - i\xi ) /\sqrt{2}$,
we obtain the desired formula for $I (g)$. Next following the
stationary phase method, we prove that the map $I$ is continuous $S(0)
\rightarrow S(0)$ with the asymptotic expansion
(\ref{equ:toeplitz-weyl}). Using a density argument, we conclude that
the Weyl symbol of $U^*_\B T_g U_\B $ is $I(g)$ for any $g$ in the
class $S(0)$.
\end{demo}

Since the series involved in~\eqref{equ:toeplitz-weyl} is that of the
exponential, it is easy to inverse formally. Hence let $g(\cdot,h)$ be
a family in $\Cinf_0(\CM^n\times [0,\epsilon_0])$ with support in $D$
that admits the following asymptotic expansion:
\begin{equation}
  g(\cdot,h) = \sum(\sfrac{-h}{2})^j \Delta^jq/j! + \O(h^\infty) .
  \label{equ:g}
\end{equation}
Then the proposition says that
\[
Q^\B = T_{g(\cdot;\epsilon,h)} + \O(h^\infty)
\]
where the remainder is in the uniform norm. Notice that, since
$h^{-1}$ is an affine function of $N$ with positive slope, asymptotics
in $h$ are equivalent to asymptotics in $N^{-1}$. In particular we can
neglect the remainder $\O(h^\infty)$. Therefore, restricting to
$\H_h^\B$ we can express $K^\B(\epsilon,h)=Q^\B_{\restr \H_h^\B}$ as
\begin{equation}
  K^\B(\epsilon,h):\H_h^\B\fleche \H_h^\B, \qquad \psi\mapsto
  \BPi^{\H_h^{\B}} g(\cdot;\epsilon,h)\psi,
  \label{equ:KB}
\end{equation}
where $\BPi^{\H_h^{\B}}$ is the orthogonal projector of
$L^2(\CM^n,e^{-\abs{z}^2/h})$ onto $\H_h^\B$.

It remains to switch to the space $\BH_N$ and see how to reduce the
symbol $g=g(\cdot;\epsilon,h)$ to the orbifold $M$.

The first step is to prove the theorem with $(V_N^*)^{-1}$ instead of
$U_N$. In other words we look for a suitable symbol $J(g)\in \Cinf(M)$
such that
\[
(V_N^*)^{-1} K^\B(\epsilon,h) V_N^{-1} = \BPi_N J(g),
\]
or, equivalently, inserting~\eqref{equ:KB},
\[
\BPi^{\H_h^{\B}} g = V_N^* \BPi_N J(g) V_N,
\]
acting on $\H_h^\B$. Using that $\BPi^{\H_h^{\B}}$ and $\BPi_N$ are
self-adjoint, this amounts to show that for any $\psi,\psi'\in\H_h^\B$
\begin{equation}
  \pscal{g\psi}{\psi'}_{\CM^n} = \pscal{J(g)V_N\psi}{V_N\psi'}_M.
  \label{equ:pscalaire}
\end{equation}

By definition
\begin{equation}
  \pscal{g\psi}{\psi'}_{\CM^n} = \int_{\CM^n} e^{-\abs{z}^2/h}
  g(z)\psi(z)\overline{\psi'(z)}\mu(z),
  \label{equ:integrale}
\end{equation}
where $\mu$ is Lebesgue's measure on $\CM^n=\RM^{2n}$. In order to
decompose this integral, we use a slicing of $\CM^n$ transversal to
$Y$ given by the flow of the harmonic oscillator $H_2$ at imaginary
times. Precisely, the map
\[
\RM \times Y \fleche \CM^n\setminus\{0\}, \quad (t,y)\mapsto
z=(e^t.y):=(e^{t\p_1}y_1,\dots,e^{t\p_n}y_n)
\]
is a diffeomorphism. Now recall that $\H_h^\B$ is spanned by
$z^\alpha$, $\pscal{\alpha}{\p}=N$. So any element $\Psi \in \H_h^\B$
satisfies
\[
\psi (u.z) = u^N\psi(z), \qquad \forall u\in\CM.
\]
Hence when $u=e^t$ one can write
\begin{equation}
  \abs{\psi(e^t.z)}^2 e^{-h^{-1}\abs{e^t.z}^2} = \abs{\psi(z)}^2
  e^{-\nu_cN\abs{z}^2} e^{-N\phy(t,z)}
  e^{-\frac{\abs{\nu}}2\abs{e^t.z}^2}
  \label{equ:phase1}
\end{equation}
with
\[
\phy(t,z) = -2t +\nu_c \sum (e^{2\p_i t}-1)\abs{z_i}^2.
\]
Since the origin is not in the support of $g$ we shall now use the
coordinates $(t,y) \in \RM \times Y$ in order to calculate the
integral~\eqref{equ:integrale}. The measure $\mu$ can be decomposed as
$$\mu=\delta(t,y) \abs{dt}\mu_Y (y),$$ where $\delta(t,y)$ is smooth and
$\mu_Y$ is the Liouville measure of $Y$. Recall that $\mu_Y$ is
$S^1$-invariant and $\pi_*\mu_Y=\mu_M$.  Since $\mu$ is also
$S^1$-invariant, the function $\delta(t,y)$ must be $S^1$-invariant as
well. Finally remark that $g$ also, as defined in~\eqref{equ:g},
inherits the $S^1$-invariance of $q$. This entails, together
with~\eqref{equ:phase1}, that
\[
\pscal{g\psi}{\psi'}_{\CM^n} = \int_Y e^{-\nu_c N\abs{y}^2}
J^Y(g)(y)\psi(y)\overline{\psi'(y)}\mu_Y(y),
\]
where the function $J^Y(g)\in\Cinf(Y)$ is defined by
\[
J^Y(g)(y) = \int_\RM e^{-N\phy(t,y)} g(t,y) e^{-\frac{\nu}2\abs{ty}^2}
\delta(t,y) \abs{dt} .
\]
Since $\phy$ is $S^1$ invariant, so is $J^Y(g)$. Therefore there
exists $J(g)\in\Cinf(M)$ such that $\pi^*J(g)= J^Y(g)$. Recall the
definition (\ref{eq:produit_scalaire}) of the scalar product of
$\BH_N$. Since $\pi_*\mu_Y=\mu_M$, we get the desired
identity~\eqref{equ:pscalaire}.

Asymptotics of $J^Y(g)$ (and hence of $J(g)$) are obtained by the
stationary phase lemma. Since $\phy(t,y)$ has a global non-degenerate
minimum at $t=0$, the expansion is localised on $Y$, as expected. For
instance, at first order, we see that if $f_0\in\Cinf(M)$ is such that
$\pi^*f_0(y)=g(0,y;\epsilon,h)+\O(h)$ then
\[
J(g)(m) = e(m)f_0(m) +\O(h)
\]
where $e\in\Cinf(M)$ is a positive function.

To complete the proof of the theorem it remains to replace
$(V_N^*)^{-1}$ by $U_N$. But since $U_N=(V_NV_N^*)^{-\frac12}V_N$ we
have
\[
U_N K^\B(\epsilon,h) U_N^* = (V_N V_N^*)^{\frac12} (V_N^*)^{-1}
K^\B(\epsilon,h) V_N^{-1} (V_N V_N^*)^{\frac12},
\]
which means
\[
U_N K^\B(\epsilon,h) U_N^* = (V_N V_N^*)^{\frac12} \left(\BPi_N
  J(g)\right) (V_N V_N^*)^{\frac12}.
\]
If we repeat our argument with $K(\epsilon,h)=
\operatorname{Id}+\O(h^\infty)$ we see that
$$(V_N^*)^{-1}( \operatorname{Id} + \O(h^\infty))V_N^{-1}$$ is a Toeplitz operator, and its principal symbol is $e$.  We cannot omit here the remainder
$\O(h^\infty)$ because we work with a compactly supported symbol for
$K$. However, one can prove that the uniform norms of $
N^{-\frac{1}{4}} V_N$ and its inverse are $\O (1)$ (cf.
~\cite{charles-reduc}, proposition 4.22) Hence
$(V_N^*)^{-1}(\O(h^\infty))V_N^{-1}=\O(h^\infty)$ and thus is a
Toeplitz operator (with asymptotically trivial symbol).
Hence $W_N:=(V_N^*)^{-1}V_N^{-1}$ is a Toeplitz operator. By the
symbolic and functional calculus of Toeplitz operator we get that
\[
U_N K^\B(\epsilon,h) U_N^* = W_N^{-\frac12} \left(\BPi_N J(g)\right)
W_N^{-\frac12}
\]
is indeed a Toeplitz operator with principal symbol
$e^{-\frac12}(ef_0)e^{-\frac12}= f_0$.
\end{demo}

\subsection{Spectral asymptotics of eigenvalue clusters}
\label{sec:asymptotics}

In this section we apply the previous result to the Birkhoff normal
form of a pseudo-differential operator $P$, in order to get spectral
asymptotics that we express in the original parameters $(E,\h)$. So we
assume that $P$ satisfies the assumptions of theorem \ref{theo:FNB}
with
$$H_2(x,\xi)= \frac{\nu_c}{2}\sum_{j=1}^n \p_j (x_j^2+\xi_j^2) .$$
For the sake of simplicity, we shall also assume that $H_0=0$~: the
sub-principal symbol of $P$ vanishes at $z_0$.  Then, according to
this theorem, the small eigenvalues of $P$ correspond to eigenvalues
of $\hat{H}_2 + K$, where $K$ commutes with $\hat{H}_2$.

Let $r$ be equal to $3$ when there is a resonance relation of the form
$\p_j=2\p_i$ or $\p_i = \p_j + \p_k $, and to $4$ otherwise. Then $K$
belongs to $\M_{r,0}$ and by lemma \ref{lemm:normes}, the norm of the
restriction of $K$ to the eigenspace $\H_{E,\h}$ of $\hat{H}_2$ is $\O
( E^ \frac{r}{2})$.

Since the distance between two consecutive eigenvalues of $\hat{H}_2$
is larger than $\hbar \nu_c$, we conclude that the bottom of the
spectrum of $P$ split into bands or clusters in this precise sense:
there exists $\hbar_0>0$ and $C>0$ such that for every $\hbar \in (0,
\h_0]$
\begin{equation}
  \Sp ( P ) \cap (-\infty, C \hbar ^{\frac{2}{r}} ) \subset \bigcup_{E
    \in \Sp ( \hat{H}_2 )} \Bigl[ E - \frac{\nu_c \hbar}{3}, E +
  \frac{\nu_c\hbar}{3} \Bigr]
  \label{equ:clusters}
\end{equation}
Furthermore for any eigenvalue $E$ for $\hat{H}_2$ smaller than $C
\hbar ^{\frac{2}{r}}$,
\begin{equation*}
  \# \Sp (P) \cap  \Bigl[ E - \frac{\nu_c \hbar}{3}, E +
  \frac{\nu_c\hbar}{3} \Bigr] = m (E, \h)
\end{equation*}
with the multiplicity $m(E, \h)=\dim \H_{E,\h}$. First result is an
estimate of the width of the band and of the distribution of the
eigenvalues in each band in the regime $\h / E \rightarrow 0$.  Denote
by
$$ E + \lambda_1 (E, \h ) \leqslant .... \leqslant E +
\lambda_{m(E,\h)} (E, \h)$$ the eigenvalues of $P$ in $\bigl[ E -
\frac{\nu_c \hbar}{3}, E + \frac{\nu_c\hbar}{3} \bigr]$.

Let $k_0$ be the principal $\h$-symbol of $K$.

\begin{theo}
  \label{theo:int1}
  For any $\h \in (0, \h_0]$ and any eigenvalue $E \leqslant C
  \h^{\frac{2}{r}}$ for $\hat{H}_2$, we have
  $$ \lambda_1 (E, \h)  = \inf_{x \in \{H_2 = E \}} \abs{k_0(x)} +
  E^{\frac{r}{2}} \O ( \h / E ), $$
  $$\lambda_{m(E,\h)} (E, \h)  = \sup_{x \in \{H_2 = E \}} \abs{k_0(x)} +
  E^{\frac{r}{2}} \O ( \h/ E)
  $$
  and for any function $g \in \Cinf ( \RM)$,
  $$ \sum_{i=1}^{m (E, \h)} g \bigl( \lambda_i (E, \h)/ E^{\frac{r}{2}} \bigr) =
  \Bigl( \frac{1}{2\pi \h}\Bigr)^{n-1} \int_{\{ H_2 = E\}} g \Bigl(
  \frac{k_0(x)}{E^{\frac{r}{2}}} \Bigr) \mu_E (x) + \O(( E/ \h)^{n-2})
  $$ where $\mu_E$ is the Liouville measure of $\{ H_2 = E \}$ and the
  $\O$'s are uniform with respect to $\h$ and $E$.
\end{theo}

To read this result, it is interesting to understand the dependence of
the leading order terms with respect to $\h$ and $E$ as well. Recall
that by remark \ref{rema:class_birkhoff} the Taylor expansion of $k_0$
is precisely the classical Birkhoff normal form of the principal
symbol of $P$. We have $k_0= k^{(r)} +\O ( |x,\xi|^{r+1})$ with
$k^{(r)} (x, \xi)$ a homogeneous polynomial of degree $r$. Then
\[
\inf_{\{H_2 = E \}} \abs{k_0} = E^{\frac{r}{2}} \inf_{\{H_2 = 1 \}}
\abs{k^{(r)}} + \O ( E^{\frac{r +1}{2}}),
\]
\[
\int_{\{ H_2 = E\}} g \Bigl( \frac{k_0(x)}{E^{\frac{r}{2}}} \Bigr)
\mu_E (x) = E^{n-1} \int_{\{ H_2 = 1\}} g (k^{(r)}(x) ) \mu_1 (x) + \O
( E^{ n - \frac{1}{2}}) .
\]

If $r =3$, $k^{(r)}$ is easily computed. Let the Taylor expansion of
the principal symbol of $P$ begin with $ H_2 + p^{(3)} + \O ( |x,
\xi|^4).$ Then $k^{(3)}$ is the average of $p^{(3)}$ with respect to
the Hamiltonian flow of $H_2$. If $r=4$ then in general the formula is
more involved~: letting $p=H_2+p^{(3)}+p^{(4)} + \O(\abs{x,\xi}^5)$ we
see from the Birkhoff construction that $k^{(4)}$ is the average of
$p^{(4)}-\frac12\{p^{(3)},a^{(3)}\}$, where $a^{(3)}$ is the
Hamiltonian whose flow performs the first averaging in the method,
\emph{ie.} $\{H_2,a^{(3)}\}=-p^{(3)}$. We do not know a simpler
formulation for this term, except, of course, when $p^{(3)}=0$.

\begin{demo}[of the theorem]
  First the result holds for the eigenvalues
  \[
  \lambda_1' (E,\h) \leqslant ...\leqslant \lambda_{ m (E, \h)}' (E,
  \h)
  \]
  of the restriction of $K$ to the eigenspace $\H_{E,\h}$.  Indeed, by
  theorem \ref{theo:toeplitz}, the numbers $\epsilon^{-r} \lambda_i' (
  \epsilon ^2, \epsilon^2 h)$ are the eigenvalues of a Toeplitz
  operator which depends smoothly on $\epsilon$ and with semiclassical
  parameter $h$. The principal symbol of this operator is the
  push-forward to $M$ of the restriction of $\epsilon^{-r} k_0 (
  \epsilon .) $ to $\{ H_2 = 1\}$ or equivalently of the restriction
  of $\epsilon^{-r} k_0 (.) $ to $\{ H_2 = \epsilon^2\}$. The
  semiclassical estimates of the smallest and largest eigenvalue and
  of the spectral density of a Toeplitz operator in terms of its
  principal symbol are basic results which were extended to the
  orbifold case in \cite{charles-reduc}. Going back to the original
  parameter $E$ and $\h$, we obtain the result for the $\lambda_i'(E,
  \h)$ and the $\O$ are uniform when $E$ and $\h$ run over an
  arbitrary bounded set.
  
  Then when $E \leqslant C \h^{\frac{2}{r}}$ and $\h \in (0, \h_0]$ we
  have for any $N$
  $$  \abs{ \lambda_i (E,\h)-\lambda_i' (E,\h)} \leqslant C_N (E^N +
  \h^N) $$ and we conclude easily. Since $\frac{1}{2} \nu_c \h
  \leqslant E \leqslant C \h^{\frac{2}{r}}$ the remainders $\O(E^N)$
  and $\O(\h^N)$ are negligible when $N$ is sufficiently large and
  disappear.
\end{demo}

\begin{rema}
  \label{rema:polyade}
  In some sense the result is still true for higher energy $E$. When
  the condition $E\leq C\h^{\frac{2}{r}}$ is violated, the bands may
  overlap, and we can not extract the eigenvalues $\lambda_i (E,\h)$
  from the whole spectrum of $P$ as we did in~\eqref{equ:clusters}.
  However we saw in the proof that the eigenvalues $\lambda_i'(E,\h)$
  satisfy the same estimates as the $\lambda_i (E, \h)$ except that
  the remainders are uniform with respect to bounded energy $E$.
  Therefore the theorem remains valid in this regime provided
  $\lambda_i (E, \h)$'s are replaced by $\lambda_i'(E,\h)$'s.  Then we
  can recover the spectrum of $P$ with theorem \ref{theo:FNB}.
\end{rema}
\begin{figure}[h]
  \centering
  \includegraphics[width=\textwidth]{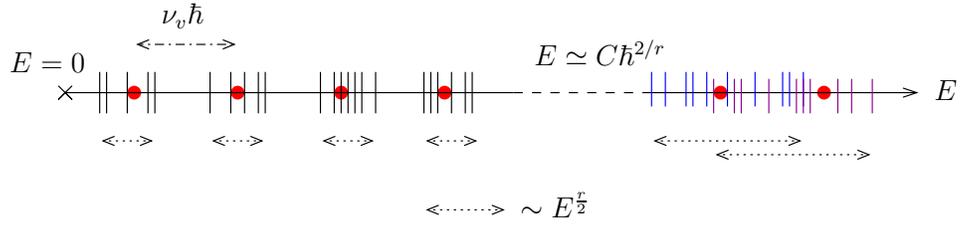}
  \caption{\small The spectrum (vertical ticks) exhibits clustering
    around the eigenvalues of the harmonic oscillator (big dots) when
    $E\leq C\h^{2/r}$.  But even when overlapping occurs, one can give
    a sense to each band if one thinks in terms of the eigenvalues
    $\lambda_i'(E,\h)$.  See remark~\ref{rema:polyade}.}
  \label{fig:clusters}
\end{figure}

The estimate of the spectral density in the previous theorem is
actually the leading order term of a full asymptotic expansion. The
description of its various pieces is involved because the reduced
phase space $M$ is not a manifold but an orbifold.  First denote by
$G$ the set of complex number $\zeta$ of modulus 1 whose fixed point
set $$ Y_\zeta : = \{ z \in Y; \; \zeta.z = z \}$$ is not empty. Here
the dot is for the Hamiltonian action generated by $H_2$ as in
(\ref{eq:action}).  A straightforward computation leads to
\begin{gather*}
  G = \{ \zeta \in \CM^*; \; \zeta^{\mathfrak{p}_i} = 1 \text{ for
    some }i \}
\end{gather*}
and
\begin{equation}
  Y_\zeta = \CM_{\zeta} \cap Y \text{ with }\CM_{\zeta} = \{ z \in
  \CM^n ; \; z_i =0 \text{ if } \zeta ^{\p_i} \neq 1 \}.
  \label{equ:Y}
\end{equation}
The $S^1$-action preserves $Y_\zeta$. Let $M_\zeta$ be the quotient
$Y_\zeta / S^1$. It is a twisted projective space which embeds into
$M$ as a symplectic suborbifold. Denote by $n(\zeta)$ its complex
dimension. Finally, let $m(\zeta)$ be the greatest common divisor of
$\{ {\mathfrak{p}}_i; \; \zeta^{\mathfrak{p}_i}=1\}$. Then with the same argument
as in the previous proof we deduce from theorem 2.3 of
\cite{charles-reduc} the

\begin{theo}
  \label{theo:int2}
  For every function $g\in \Cinf( \RM)$,
  \begin{gather*}
    \sum_{i=1}^{m (E, \h)} g \bigl( \lambda_i (E, \h)/ E^{\frac{r}{2}}
    \bigr) = \sum_{ \zeta \in G} \Bigl( \frac{E}{2\pi \h}
    \Bigr)^{n(\zeta )} \zeta^{-N} \sum_{\ell=0}^{\infty} \Bigl(
    \frac{\h}{E}\Bigr)^{\ell} I_\ell( E^{\frac{1}{2}}, \zeta ) + \O
    \Bigl( \Bigl(\frac{\h}{E}\Bigr)^{\infty}\Bigr)
  \end{gather*}
  where the remainder is uniform when $\h$ and $E$ runs over $(0,
  \h_0)$ and $\Sp (\hat{H}_2) \cap (-\infty, C \h^{\frac{2}{r}})$
  respectively, and $N$ is defined by $E=\h(\frac{\nu}2+N)$.
  Furthermore the coefficients $ I_\ell(\epsilon,\zeta)$ are $\Cinf$
  function of $\epsilon$ and
$$ I_0( \epsilon, \zeta ) = \frac{1}{m(\zeta)} \Biggl( \prod_{i; \;
  \zeta^{{\mathfrak{p}}_i} \neq 1} (1 - \zeta^{\mathfrak{p}_i} )^{-1} \Biggr)
\int_ {M_{\zeta}} g(\epsilon^{-r} k_0(\epsilon x ) ) \mu_{\zeta} (x),
$$ where $ \mu_ {\zeta}$ is the Liouville measure of $M_{\zeta}$.
\end{theo}

\subsection{A convex polytope and a trace formula}

As a consequence of the preceding theorems~\ref{theo:int1}
and~\ref{theo:int2} we obtain an interesting formula expressing the
asymptotics of a combinatoric sum over integral points of a convex
polytope, when the polytope undergoes some rational scaling.

We first state the result without any reference to any operator,
recalling only the following notation:

Let $n$ be a positive integer and let $\p_1,\dots,\p_n$ be coprime
positive integers. Let $\p=(\p_1,\dots,\p_n)$.
Again let
\begin{gather*}
  G = \{ \zeta \in \CM^*; \; \zeta^{\p_i} = 1 \text{ for some }i \}
\end{gather*}
and for each $\zeta\in G$ introduce the index set
$\mathbf{i}_\zeta=\{i; \zeta^{\p_i}=1\}$, whose cardinality is
$n(\zeta)+1$.
Finally recall that $m(\zeta)=\gcd \{\p_i; \; i\in\mathbf{i}_\zeta\}$.

For $N\in\NM$ and $\alpha\in\NM^n$, we define the convex polytope
(actually a simplex of dimension $n-1$)

\[
\P(\alpha,N) = (\RM^+)^n\cap \{(x_1,\dots,x_n); \quad
\pscal{x}{\p}=N-\pscal{\alpha}{\p}\}.
\]
\begin{rema}
  Notice that $\P(\alpha,N)$ is neither integral or Delzant (in the
  terminology of \cite{guillemin-moment-book}). However if $p$ is the
  least common multiple of $\p_1,\dots,\p_n$ then $\P(p\alpha,pN)$ has
  vertices with integral coordinates, but it still not Delzant in
  general. This reflects the fact that $\P(\alpha,N)$ is a moment
  polytope for a symplectic orbifold which --- except for
  $\p_1=\dots=\p_n=1$ --- is not a manifold.
\end{rema}
When $\alpha=(\alpha_1,\dots,\alpha_n)\in\NM^n$ we use the notation
$\alpha!=(\alpha_1)!\cdots(\alpha_n)!$.
\begin{theo}
  \label{theo:polytope}
  For any $\alpha\in\NM^n$ we have the asymptotic formula
  \begin{gather}
    \label{equ:polytope}
    \frac{1}{N^{\abs{\alpha}}}\!\!\!\!\!\!\!\
    \sum_{\gamma\in\P(\alpha,N)\cap\ZM^n}
    \frac{(\gamma+\alpha)!}{\gamma!} = \sum_{\zeta\in G} \zeta^{-N}
    N^{n(\zeta)} \sum_{\ell=0}^\infty N^{-\ell} a_{\ell}(\alpha,\zeta)
    + \O(N^{-\infty})~.
  \end{gather}
  Moreover
  \[ a_{0}(\alpha,\zeta) = \frac{\displaystyle
    \prod_{i\not\in\mathbf{i}_\zeta} (1 - \zeta^{p_i}
    )^{-1}}{\displaystyle m(\zeta)\prod_{i\in\mathbf{i}_\zeta}
    \p_i^{\alpha_i+1}} \frac{\displaystyle
    \prod_{i\in\mathbf{i}_\zeta} \Gamma(\alpha_i+1)}
  {\Gamma\Bigl(\sum_{i\in\mathbf{i}_\zeta}(\alpha_i+1)\Bigr)}~.
  \]
\end{theo}
\begin{rema}
  As we shall see below, the $a_0(\alpha,\zeta)$'s are actually ``Weyl
  terms'', in the sense that they are the result of some integrals
  over different faces of the polytope (or, equivalently, they are
  phase space integrals for some sub-orbifolds).
\end{rema}
\begin{rema}
  If we choose $\alpha=0$, the left-hand-side of is just the number of
  integral points of the polytope $\P(0,N)$, and the formula becomes:
  \[
  \#(\P(0,N)\cap\ZM^n) = \sum_{\zeta\in G} \zeta^{-N} N^{ n(\zeta)}
  \sum_{\ell=0}^\infty N^{-\ell} b_{\ell}(\zeta) + \O(N^{-\infty})
  \]
  where
  \[
  b_0(\zeta) = a_{0}(0,\zeta) = \frac{\displaystyle
    \prod_{i\not\in\mathbf{i}_\zeta} (1 - \zeta^{p_i}
    )^{-1}}{\displaystyle m(\zeta) (\prod_{i\in\mathbf{i}_\zeta} \p_i)
    \Gamma(n(\zeta)+1)}
  \]
  It is known that each factor of $\zeta^{-N}$ is actually a
  polynomial in $N$ given by a Riemann-Roch type formula.
  See~\cite{meinrenken-RR} for details on this issue. The leading
  term, obtained with $\zeta=1$ is
  \[
  \frac{N^{n-1}}{(\p_1\cdots\p_n)(n-1)!}
  \]
  The other coefficients can also be obtained as the coefficients of
  the generating function $g(X)=\prod_{i=1}^n(1-X^{\p_i})^{-1}$;
  see~\cite{sado-zhilin-polyades}.
\end{rema}
\begin{demo}[of the theorem]
  Let us work in the Bargmann representation. The harmonic oscillator
  we consider is
  \[
  \hat{H}_2(\h) = \h\sum_{j=1}^n (\p_j\deriv{}{z_j}+\frac12).
  \]
  If $\alpha$ and $\beta$ are multi-indices in $\NM^n$, the
  differential operator $z^\alpha(\h\deriv{}{z})^\beta$ commutes with
  $\hat{H}_2(\h)$ if and only if $\pscal{\alpha-\beta}{\p}=0$. Thus
  let us consider the symmetric differential operators
  \[
  K_\alpha(\h):=z^\alpha \bigl( \h\deriv{}{z} \bigr)^\alpha.
  \]
  We shall compute the trace of the restriction of $K_\alpha$ to the
  eigenspace $\H_{E,\h}$ in two different ways. The first way is just
  to do an explicit computation in a basis of $\H_{E,\h}$. The second
  way is to use remark~\ref{rema:polyade} in order apply
  theorem~\ref{theo:int2} with $g=\textup{Id}$ and $E=1$.

  \paragraph{1. ---} A basis of $\H_{E,\h}$ (or, more exactly, of the
  space $\H_{E,\h}^\B$ in the Bargmann representation) is given by the
  monomials
  \[
  \{z^\gamma;\quad \gamma\geq 0 \textup{ and } \pscal{\gamma}{\p}=N\}
  \]
  where the integer $N$ is defined by the equation
  $E=\h(N+\abs{\p}/2)$, and the inequality $\gamma\geq 0$ stands for
  $\gamma_j\geq 0, \forall j$.  It is straightforward to check that
  \[
  K_\alpha(\h)(z^\gamma)= \begin{cases}
      \h^{\abs{\alpha}}\frac{\gamma!}{(\gamma-\alpha)!}z^{\gamma}
      & \text{ if
      } \gamma\geq\alpha\\
      0 & \text{ otherwise.}
    \end{cases}
  \]
  Hence
  \[
  \operatorname{Tr} (K_{\alpha}(\h)_{\restr \H_{E,\h}}) =
  \h^{\abs{\alpha}}\sum_{\gamma\geq \alpha \atop \pscal{\gamma}{\p}=N}
  \frac{\gamma!}{(\gamma-\alpha)!} = \h^{\abs{\alpha}}\!\!\!\!\!\!\!\!
  \sum_{\gamma\in\P(\alpha,N)\cap\ZM^n}
  \frac{(\gamma+\alpha)!}{\gamma!}.
  \]
  \paragraph{2. ---} Formula~(\ref{equ:polytope}) is now a simple
  transcription of theorem~\ref{theo:int2}, using $E=1$ and
  $\h=1/(N+\abs{\p}/2)$, and rewriting the coefficients in order to
  transform the expansion in powers of $(\frac{\h}{E})$ into an
  expansion in powers of $N^{-1}$. This of course does not change the
  formulas for the leading coefficients.

  Let us compute these leading coefficients. For this we use some
  homogeneity property in the variable $E$, so it is best not to set
  $E=1$ for the moment.

  We introduce the standard action-angle coordinates for the harmonic
  oscillator in $\RM^{2n}$. Let
  $I_j=\frac{1}{2}(x_j^2+\xi_i^2)=\abs{z_j}^2$, and let the angles
  $\theta_j$ be defined by
  \[
  \begin{cases}
    x_j = \sqrt{2I_j}\cos\theta_j\\
    \xi_j = -\sqrt{2I_j}\sin\theta_j~.
  \end{cases}
  \]
  Since $d\xi_j\wedge dx_j = dI_j\wedge d\theta_j$ and the singular
  set of these action-angle coordinates is of codimension 2, we can
  use the corresponding symplectic measure on $(\RM^+)^n\times
  (\RM/2\pi\ZM)^n$ : $\abs{dI_1\wedge\cdots\wedge dI_n\wedge
    d\theta_1\wedge\cdots \wedge d\theta_n}$ as a replacement for the
  symplectic measure $\mu$ on $\RM^{2n}$. Let $\abs{dE}$ denote the
  pull-back by $H_2=\sum_j\p_j I_j$ of the Lebesgue measure on $\RM$.
  Since the flow of $H_2$ is $2\pi$-periodic, the Liouville measure
  $\mu_E$ on $\{H_2=E\}$ is by definition the quotient of $\mu$ by
  $2\pi\abs{dE}$~:
  \[
  \mu = 2\pi \mu_E \otimes\abs{dE}.
  \]
  The principal symbol of $K_\alpha(\h)$ is
  $k_{0,\alpha}=z^\alpha\bar{z}^\alpha=\prod I_j^{\alpha_j}$, for
  which we use the notation $I^\alpha$. Let
  $S_\alpha(E)=\int_{\{H_2=E\}} I^\alpha \mu_E$.

  Since $dE\wedge dI_2\wedge\cdots\wedge dI_n=\p_1
  dI_1\wedge\cdots\wedge dI_n$ we see that
  \[
  \mu_E={(2\pi\p_1)}^{-1}\abs{dI_2\wedge\cdots\wedge dI_n\wedge
    d\theta_1\wedge\cdots \wedge d\theta_n}
  \]
  and
  \[
  S_\alpha(E)=(2\pi)^{n-1}\p_1^{-1}\int_{\p_1I_1+\cdots \p_n I_n=E}
  I^\alpha \abs{dI_2\cdots dI_n} = E^{\abs{\alpha}+n-1}S_\alpha(1).
  \]
  We conclude by adapting the standard trick used to calculate the
  surface of the unit sphere, namely let
  \[
  J=\int_0^\infty e^{-E}S_\alpha(E)\abs{dE}.
  \]
  On the one hand,
  \[
  J= S_\alpha(1)\int_0^\infty e^{-E}E^{\abs{\alpha}+n-1}\abs{dE} =
  S_\alpha(1)\Gamma(\abs{\alpha}+n);
  \]
  on the other hand
  \begin{equation}
    \begin{split}
      J &= \int_{(\RM^+)^n} e^{-H_2(I)}I_1^{\alpha_1}\cdots
      I_N^{\alpha_n}\abs{dI_1\cdots dI_n}\\
      &= \left(\int_0^\infty
        e^{-\p_1I_1}\abs{dI_1}\right)\cdots\left(\int_0^\infty
        e^{-\p_nI_n}\abs{dI_n}\right).
    \end{split}
  \end{equation}
  Hence
  \[
  S_\alpha(1)=\frac{\Gamma(\alpha_1+1) \cdots
    \Gamma(\alpha_n+1)}{\p_1^{\alpha_1+1}\cdots\p_n^{\alpha_n+1}
    \Gamma(\abs{\alpha} + n)}.
  \]
  This gives the result for $a_0(\alpha,1)$. For a general $\zeta$ the
  calculation of the integral over $M_\zeta$ is exactly the same as
  the integral $S_\alpha(1)$ provided we keep only indices
  $i\in\mathbf{i}_\zeta$, as follows from the definition of $Y_\zeta$
  in~\eqref{equ:Y}. This finishes the proof.
\end{demo}

\begin{rema}
  Instead of $K_\alpha(\h)$, we could have considered the more general
  operators commuting with $\hat{H}_2(\h)$:
  $z^\alpha(\h\deriv{}{z})^\beta$, with $\pscal{\alpha-\beta}{\p}=0$.
  But using the basis $z^\gamma$ as before it is easy to see that the
  trace of such operators always vanishes as soon as
  $\alpha\neq\beta$. As a consequence, this shows that for any
  $\zeta\in G$,
  \[
  \int_{M_\zeta} z^\alpha{\bar{z}}^\beta \mu_\zeta = 0 \qquad \text{
    provided } \pscal{\alpha-\beta}{\p}=0,\quad \alpha\neq\beta.
  \]
\end{rema}

\vspace{1cm}

\bibliographystyle{plain}\bibliography{bibli}
\end{document}